\documentclass{article}

\usepackage{PRIMEarxiv}

\usepackage[utf8]{inputenc} 
\usepackage[T1]{fontenc}    
\usepackage{hyperref}       
\usepackage{url}            
\usepackage{booktabs}       
\usepackage{amsfonts}       
\usepackage{nicefrac}       
\usepackage{microtype}      
\usepackage{lipsum}
\usepackage{fancyhdr}       
\usepackage{graphicx}       
\graphicspath{{media/}}     

\usepackage{amsmath}
\usepackage{bm}
\usepackage{xcolor}
\usepackage{subfig}

\def\rbold{\mathbf{r}}
\def\nbold{\mathbf{n}}
\def\Gbold{\mathbf{G}} 
 
\def\Hbold{\mathbf{H}} 
\def\hatHbold{\mathbf{\hat H}} 
\def\hatGbold{\mathbf{\hat G}}
\def\hatqbold{\mathbf{\hat q}} 

\def\phibold{\bm{\phi}}

\def\hatphibold{\bm{\hat\phi}}

\title{Boundary element formulation of the Mild-Slope Equation for harmonic water waves propagating over unidirectional variable bathymetries
\thanks{\textit{\underline{Citation}}: 
\textbf{Antonio Cerrato, José A. González, Luis Rodríguez-Tembleque. Boundary element formulation of the Mild-Slope Equation for harmonic water waves propagating over unidirectional variable bathymetries. Engineering Analysis with Boundary Elements, Volume 62, January 2016, Pages 22-34 DOI:10.1016/j.enganabound.2015.09.006.}} 
}

\author{
  Antonio Cerrato, José A. González, Luis Rodríguez-Tembleque \\
  Escuela Ténica Superior de Ingeniería \\
  Universidad de Sevilla \\
  Camino de los Descubrimientos s/n, 41092 Sevilla, Spain\\
  \texttt{\{antoniocerrato, japerez, luisroteso\}@us.es} \\
}

\begin{document}
\maketitle

\begin{abstract}
This paper presents a boundary element formulation for the solution of the Mild-Slope equation in wave propagation problems with variable water depth in one direction. Based on the Green's function approximation proposed by Belibassakis \cite{Belibassakis2000}, a complete fundamental-solution kernel is developed and combined with a boundary element scheme for the solution of water wave propagation problems in closed and open domains where the bathymetry changes arbitrarily and smoothly in a preferential direction. The ability of the proposed formulation to accurately represent wave phenomena like refraction, reflection, diffraction and shoaling, is demonstrated with the solution of some example problems, in which arbitrary geometries and variable seabed profiles with slopes up to 1:3 are considered. The obtained results are also compared with theoretical solutions, showing an excellent agreement that demonstrates its potential.
\end{abstract}

\keywords{Wave propagation \and Mild-slope equation \and Helmholtz equation \and Boundary element method}

\section{Introduction}

Wave propagation in variable-depth waters is a problem of significant importance in coastal engineering with applications in the design and maintenance of harbors, coastal defense works and hydrodynamic and sediment transportation studies. It is well known that the transmission of linear waves in intermediate and deep waters can be reproduced by the elliptical Mild-Slope Equation (MSE) which was derived by Berkhoff \cite{Berkhoff1972} in the early 70s. The MSE considers simultaneously the effects of diffraction, refraction, reflection and shoaling of linear water surface waves and it is formally valid for slowly varying sea bed slopes, i.e. $\nabla h << kh$, being $h$ the water depth and $k$ the wave number. The validity of the MSE has been evaluated by Tsay and Liu \cite{Tsay1983} demonstrating that it produces accurate results for bottom slopes up to 1:1 when waves are propagating perpendicularly to the bathymetry contour lines. Nevertheless, Booij \cite{Booij1983} verified that, for general directions of wave propagation, the MSE is able to provide acceptable accuracy for bottom profiles with slopes up to 1:3, enough for practical applications.

Some extensions of the MSE have been proposed in subsequent works. For example, a time-dependent extension of the MSE was derived by Kirby \cite{Kirby1986-JFM} for the case of waves propagating over ripple beds.  Also, an Extended Mild-Slope Equation (EMSE) was proposed by Massel \cite{Massel1993} that includes higher-order terms, providing a better accuracy for more complicated bathymetries. Energy dissipation effects, such as wave breaking and bottom friction, were included in \cite{Maa2002}. Chamberlain and Porter \cite{Chamberlain1995} suggested a Modified Mild-Slope Equation (MMSE), later improved by Porter and Staziker \cite{Porter1995}, which retains the second order terms discarded by Berkhoff in the formulation of the MSE. On the other hand, Suh et al. \cite{Suh1997} derived a time-dependent equation for wave propagation on rapidly varying topography and Chandrasekera et al. \cite{Chandrasekera1997} included terms for relatively steep and rapidly undulating bathymetries. Later, Lee et al. \cite{Lee1998} presented an hyperbolic MSE for rapidly varying topography, followed by the works of Copeland \cite{Copeland1985} and Massel \cite{Massel1993} in the same direction. Finally, the recent works of Hsu et al. \cite{Hsu2001} and Li et al. \cite{Li1994,Hsu2000} considered higher-order bottom effect terms to account for a rapidly varying topography and wave energy dissipation in the surf zone. Basically, all these formulations introduce higher-order terms in the MSE due to the bottom effects, usually proportional to the square of the bottom slope or the bottom curvature.

In general, the MSE represents the basic framework for the simulation of surface wave transmission problems in variable water depths and different numerical solution procedures have been proposed in the literature since the pioneering work of Berkhoff \cite{Berkhoff1972}.

Traditionally, the MSE has been solved using the Finite Element Method (FEM) \cite{Berkhoff1976} and the Finite Difference Method (FDM), where we can include the works of Li and Anastasiou \cite{Li1992}, Panchang and Pearce \cite{Panchang1991}. Nevertheless, finite difference schemes and the finite element method present a common deficiency; open and partially reflecting boundary conditions are difficult to represent. These deficiencies have been studied by many authors, like Chen et al. \cite{Chen1974,Chen1986} using hybrid FEM formulations, together with the initial proposals of Berkhoff \cite{Berkhoff1976} and Tsay et al. \cite{Tsay1983,Tsay1989} including bottom friction effects. For the closing boundary conditions, Bettess and Zienkiewicz \cite{Bettess1978} and Lau and Ji \cite{Lau1989} used infinite elements in the outer regions. Dirichlet to Neumann (DtN) boundary conditions were proposed by Givoli et al. \cite{Givoli1990,Keller1989,Givoli1991} as an analytical procedure to reproduce exact non-reflecting boundary conditions in some particular cases. This idea, was followed by Bonet \cite{Bonet2013} to derive the discrete non-local (DNL) boundary condition. More rudimentary iterative methods have also been proposed to define absorbing boundary conditions; see Beltrami et al. \cite{Beltrami2001}, Steward and Panchang \cite{Steward2001}, Chen \cite{Chen2002} or Liu et al. \cite{Liu2008}, among others. It is important to mention that a boundary element formulation of the MSE for open domains and variable bathymetry, would be able to palliate the drawbacks of FEM, providing a better approximation for the simulation of absorbing boundaries.

The MSE problem has also been solved using the Boundary Element Method (BEM). Boundary element techniques prove to be very accurate in wave refraction-diffraction problems with open domains, presenting the additional benefit that the radiation condition to infinity is automatically satisfied. In order to improve the solution of the FEM schemes, Hauguel \cite{Hauguel1978} and Shaw and Falby \cite{Shaw1978} first coupled FEM and BEM. Hamanaka \cite{Hamanaka1997} proposed a genuine BEM based boundary condition for open, partial reflection and incident-absorbing boundaries. At the same time, Isaacson and Qu \cite{Isaacson1990} introduced a boundary integral formulation to reproduce the wave field in harbors with partial reflecting boundaries and Lee et al. \cite{Lee2002,Lee2009} included the effect of incoming random waves. The Dual Reciprocity Boundary Element Method (DRBEM) has been used to model wave run-ups by Zhu \cite{Zhu1993}. Later, this technique was extended to model internal regions with variable depth surrounded by exterior regions with constant bathymetry \cite{Liu2003,Zhu2000,Zhu2009,Hsiao2009}. More recently, Naserizabeh et al. \cite{Naserizadeh2011} proposed a coupled BEM-FDM formulation to solve the MSE in unbounded problems.  

In this context, this paper presents a BEM formulation for the MSE in wave propagation problems with variable water depth in one direction. Based on the Green's function approximation proposed by Belibassakis \cite{Belibassakis2000}, a complete fundamental-solution kernel is developed and combined with a boundary element scheme for the solution of water wave propagation problems in closed and open domains where the bathymetry changes arbitrarily and smoothly in a preferential direction. This particular case is of high practical interest, because the bathymetric lines can usually be considered straight and parallel to the coast-line. A BEM formulation of the MSE for variable bathymetry not only extends the range of applications of the BEM for the solution of coastal engineering problems but also, combined with the FEM and used as a matching condition, offers the possibility of modeling very accurately the radiation condition to deeper waters.

The paper is organized as follows. Section \ref{SECTION: Formulation MSE} first reviews the formulation of the MSE. In Section \ref{sec-FS}, the fundamental solution of the MSE for variable water depth is approximated in the frequency domain. The mathematical and numerical principles of the BEM for wave scattering problems are covered in Section \ref{sec-BEM}. Section \ref{SECTION:NUMERICAL EXAMPLES} is dedicated to the validation of the proposed BEM formulation through the solution of wave propagation problems in variable water depth. Finally, Section \ref{SECTION: CONCLUSION} closes with the conclusions.

\section{The Mild-Slope Equation}
\label{SECTION: Formulation MSE}
The classical MSE \cite{Berkhoff1972,Berkhoff1976} is obtained from the linear wave theory using a Cartesian coordinate system with the $(x,y)$-plane located on the quiescent water surface and the $z$ direction pointing upwards. Under the assumption of potential flow and integrating the velocity potential in the vertical direction with appropriated boundary conditions, the velocity potential of the water surface can be represented in the form:
\begin{equation}
\Phi(x,y,t)=\phi(x,y) e^{-i\omega t},
\end{equation}
being $i$ the imaginary unit and $t$ the time variable. This potential has to satisfy the homogeneous MSE, that may be written as: 
\begin{equation}
\nabla \cdot (cc_{g}\nabla \phi) + \omega \dfrac{c_{g}}{c} \phi = 0,
\label{eqn-MSE}
\end{equation}
where $\nabla=(\partial_{x}, \partial_{y})$ is the gradient operator, $c$ is the wave velocity and $c_{g}$ the group velocity. The water depth function $h(x,y)$, wave number $k$ and angular frequency $\omega$ of the waves are related by the dispersion equation:
\begin{equation}
\label{eqn:dispersionEquation}
\omega^{2}=gk \tanh(kh),
\end{equation}
being $g$ the gravitational acceleration ($g=9.81m/s^{2}$). This means that, for a fixed frequency and variable bathymetry, the wave number $k(x,y)$ is a function of the local water depth.

The MSE can be simplified introducing the following change of variable due to Bergmann \cite{Bergmann1946}:
\begin{equation}
\label{eqn:ChangeVariable}
\phi=\dfrac{1}{\sqrt{cc_{g}}}\hat{\phi},
\end{equation}
a relation that transforms \eqref{eqn-MSE} into a Helmholtz equation: 
\begin{equation}
\label{eqn:Helmholtz}
\nabla^{2} \hat{\phi} + \hat{k}^{2}  \hat{\phi} = 0, 
\end{equation}
with a modified wave number $\hat{k}(x,y)$ given by:
\begin{equation}
\label{eqn:k con gorro}
\hat{k}^{2}(x,y)= k^{2} - \frac{\nabla^{2}\sqrt{cc_{g}}}{\sqrt{cc_{g}}},
\end{equation}
that is a known function of the wave characteristics and the local water depth.

Note that this approach is also valid for treating the same problem in the framework of the MMSE. Simply by modifying the expression of the wave number \eqref{eqn:k con gorro}, including additional effects associated with higher-order contributions of bottom slope and curvature, we obtain the MMSE model that extends the applicability of the MSE.

\section{Fundamental solution for variable wave number}
\label{sec-FS}

\begin{figure}
\centering
\def\svgwidth{0.8\columnwidth}
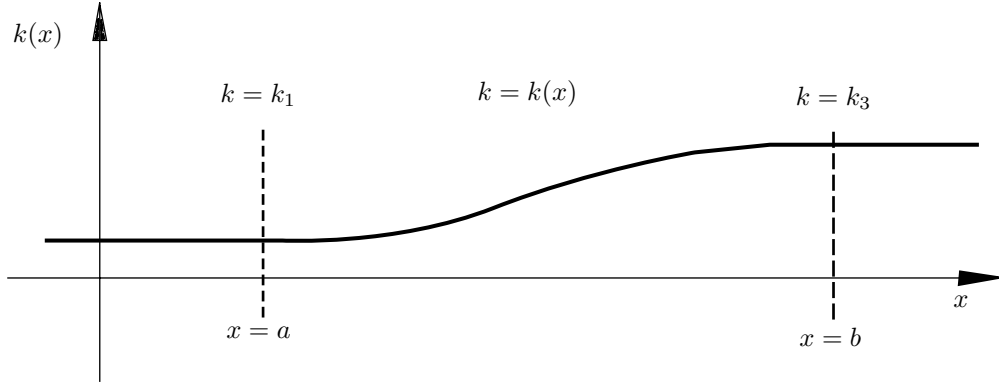
\caption{Wave number variation in the x-direction for a fixed wave frequency due to a monotonically decreasing water depth profile $h(x)$. Wave number is higher where water depth is lower as dictated by the dispersion relation}
\label{fig:Wavenumber-variation}
\end{figure}

Based on the Green's function of Belibassakis \cite{Belibassakis2000}, in this section we develop a fundamental solution of the Helmholtz problem \eqref{eqn:Helmholtz} for the particular case of an unidirectional variable bathymetry like the one described in Figure \ref{fig:Wavenumber-variation}. Taking the $x$-axis in the same direction than the variation of the water-depth $h=h(x)$, a modified wave-number $\hat{k}=\hat{k}(x)$ is obtained for a fixed wave frequency after applying relation \eqref{eqn:k con gorro}. The Green's function of the MSE equation $\psi=\psi(\rbold,\rbold_{o};\hat{k})$ is formulated as the solution of the following problem: 
\begin{equation}
\label{eqn:Green}
\nabla^{2}\psi+\hat{k}(x)^{2}\psi + \delta(\rbold-\rbold_o) = 0 \quad \text{in } \mathbb{R}^2
\end{equation}
being $\rbold_{o}=(x_{o},y_{o})$ the location of the source load and $\rbold=(x,y)$ the observed point where the velocity potential is going to be evaluated. An additional condition is that the velocity potential $\psi(\rbold)$ should satisfy the Sommerfeld's radiation condition at infinity.

To solve this problem, we apply the Fourier transform to the velocity potential in the $y$-direction, where the modified wave-number $\hat{k}(x)$ is constant, to operate with a transformed velocity potential $\varPsi=\mathcal{F}(\psi)$ that is defined as:
\begin{equation}
\label{eqn:FT-psi}
\varPsi (x,x_o;\xi) = \int_{-\infty}^{\infty} \psi(\rbold,\rbold_o; \hat{k}) e^{-i y \xi} \: dy,
\end{equation}
obtained by decomposition of the potential into its frequencies, represented by the Fourier parameter $\xi$. Introducing this transformation into problem \eqref{eqn:Green}, a family one-dimensional wave equations is obtained:
\begin{equation}
\label{eqn:1Dvarphi}
\varPsi_{,xx} + \kappa^{2}(x) \varPsi + \delta(x - x_o) = 0 \quad \text{in } \mathbb{R}
\end{equation}
with a transformed wave-number $\kappa^{2}(x)=\hat{k}^{2}(x)-\xi^{2}$ that is now function of the space variable and the Fourier parameter. Note also that the transformed problem depends on the square of this Fourier parameter, so our transformed velocity potential should be symmetric with respect to $\xi$, i.e.:
\begin{equation}
\label{eqn:symetry-xi}
\varPsi(x,x_o;\xi)=\varPsi(x,x_o;-\xi),
\end{equation}
an important property that will be used later. It is also known that the analytical solution of \eqref{eqn:1Dvarphi} for a constant wave number $\kappa$ is:
\begin{equation}
\label{eqn:FS-1Dsolkcte}
\varPsi(x,x_o;\xi) = \frac{i}{2 \kappa} e^{i \kappa |x-x_o|}
\end{equation}
and that, for large values of $\xi$, it is also possible to assume that \eqref{eqn:FS-1Dsolkcte} is a good approximation of the solution of problem \eqref{eqn:1Dvarphi} for a smooth function $\kappa(x)$. Hence for ${\xi} >> \hat{k}$ we find that $\kappa \rightarrow  i\xi$ and the transformed velocity potential decays exponentially in the form:
\begin{equation}
\label{eqn:FS-exponentially}
\varPsi(x,x_o;\xi) \simeq \frac{1}{2\xi}e^{-\xi |x- x_o|} \quad \text{for } {\xi} \rightarrow \infty,
\end{equation}
expression that defines the asymptotic behavior of the transformed potential for large values of the Fourier parameter.

Next, we observe that the one-dimensional infinite domain where the transformed problem \eqref{eqn:1Dvarphi} is defined, can be divided into three different regions, as depicted in Figure \ref{fig:Wavenumber-variation}. In the first semi-infinite interval, $x \in (-\infty,a]$, the modified wave number is considered constant $\hat k = \hat k_1$; then a second finite interval $x\in[a,b]$ where the modified wave number $\hat k=\hat k(x)$ is variable, changing monotonically from $\hat k_1$ to $\hat k_3$, and finally another semi-infinite region $x \in [b,\infty)$, where the wave number remains constant $\hat k=\hat k_3$. By performing this division of space, problem \eqref{eqn:1Dvarphi} can be reduced to a BVP defined in a finite interval $[a,b]$, with appropriated matching conditions at the boundaries, written in the following way:
\begin{align}
\label{eqn:1Dvarphi in 2}
& \varPsi_{,xx} + \kappa^{2}(x) \varPsi + \delta(x - x_o) = 0 \quad \text{in } x \in [a,b] \\
\label{eqn:bc-varphi-x}
& \left\{\begin{array}{lll}
\varPsi_{,x} + i \alpha(\xi) \varPsi  =  0 & \text{in} & x=a \\
\varPsi_{,x} - i \beta(\xi) \varPsi  =  0 & \text{in} & x=b \\
\end{array}\right.
\end{align}
being $\alpha(\xi)=(\hat{k}_1^2-\xi^2)^\frac{1}{2}$ and $\beta(\xi)=(\hat{k}_3^2-\xi^2)^\frac{1}{2}$ the parameters of the two Sommerfeld's radiation boundary conditions used to close the domain.

This one-dimensional wave transmission problem can now be solved numerically for any given value of the Fourier parameter $\xi$, providing an approximation of the transformed velocity potential in the finite interval $[a,b]$. The transformed velocity potential in the semi-infinite domains $x \in (-\infty, a)$ and $x \in (b,\infty)$, can then be substituted by the analytical solution of the equivalent one-dimensional Helmholtz problem for constant wave-number:
\begin{align}
& \varPsi(x,x_o;\xi)=\varPsi(a,x_o;\xi) \: e^{-i \alpha(\xi) \vert a-x \vert} \quad \text{in } x \in (-\infty,a], \\
& \varPsi(x,x_o;\xi)=\varPsi(b,x_o;\xi) \: e^{i \beta(\xi) \vert b-x \vert} \quad \text{in } x \in [b,\infty),
\end{align}
establishing this way the continuity of the solution in the complete domain.

After solving for the transformed velocity potential, the original variables can finally be recovered via inverse Fourier transform, $\psi=\mathcal{F}^{-1}(\varPsi)$, defined as:
\begin{equation}
\label{eqn:IFT-phi}
\psi(\rbold,\rbold_o; \hat{k}) = \dfrac{1}{2\pi} \int_{-\infty}^{\infty} \varPsi (x,x_o;\xi) e^{ i y \xi} \: d\xi,
\end{equation}
for the velocity potential and
\begin{align}
\label{eqn:x-derivative-psi}
\psi_{,x}(\rbold,\rbold_o; \hat{k}) &=
\dfrac{1}{2\pi} \int_{-\infty}^{\infty} \varPsi_{,x} (x,x_o;\xi) e^{ i y \xi} \: d\xi, 
\\
\label{eqn:y-derivative-psi}
\psi_{,y}(\rbold,\rbold_o; \hat{k}) &=
\dfrac{1}{2\pi} \int_{-\infty}^{\infty} i\xi  \varPsi (x,x_o;\xi) e^{ i y \xi} \: d\xi,
\end{align}
for its spatial derivatives.

\subsection{Inverse transform in the complex plane}
\label{subsec:complex xi}
%
\begin{figure}
\centering
\def\svgwidth{0.9\columnwidth}
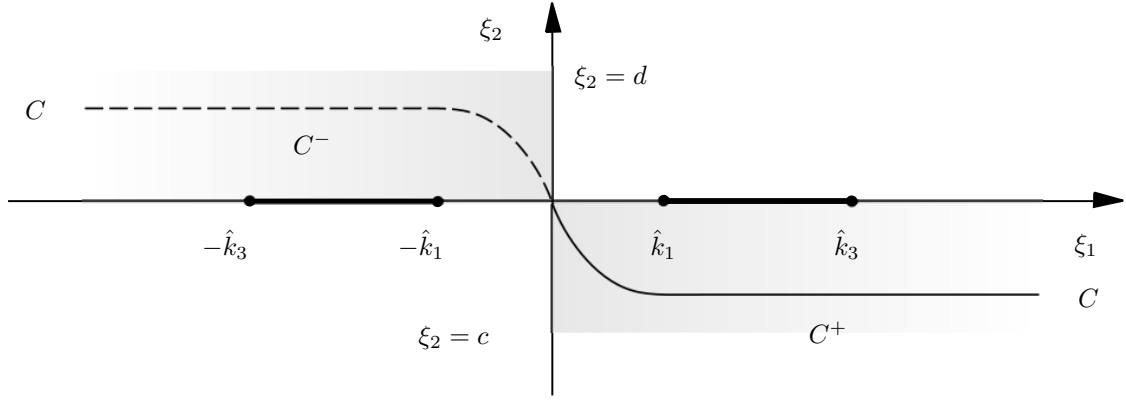
\caption{Antisymmetric integration path $C$ in the complex plane avoiding the $\pm\hat{k}_1$ and $\pm\hat{k}_3$ roots lying on the real axis}
\label{fig:IntegrationPath1}
\end{figure}
%

To compute the inverse Fourier transforms (IFT) given by equations (\ref{eqn:IFT-phi}-\ref{eqn:y-derivative-psi}), Belibassakis et al. \cite{Belibassakis2000,Belibassakis2004} propose a truncation of the infinite interval of integration and to use a discrete Fast Fourier Transform (FFT) algorithm that approximates the integral on a finite subinterval by sampling a numerical discretization of the integrand. An undesirable consequence of this approximation is that undersampling in the $\xi$-domain can cause aliasing effects in the physical $y$ direction.

In order to reduce aliasing, it is possible to extend the integration to the complex plane ($\xi=\xi_{1}+i\xi_{2}$). The selected integration path must be antisymmetric, as Equation \eqref{eqn:symetry-xi} requires, so the alternative integration path $C$ represented in Figure \ref{fig:IntegrationPath1} is used to evaluate these inverse transforms. 

Now, in the complex plane, the antisymmetry of the path $C$ is used to rewrite the IFT of the transformed velocity potential \eqref{eqn:IFT-phi} and its derivatives (\ref{eqn:x-derivative-psi}-\ref{eqn:y-derivative-psi}) in the following one-sided way:
\begin{align}
\label{eqn:IFT-phi-symmetric}
\psi(\rbold,\rbold_o; \hat{k}) &=
\dfrac{1}{\pi} \int_{\xi \in C^{+}} \varPsi (x,x_o;\xi) \cos{(y \xi)} \: d\xi, 
\\
\label{eqn:IFT-phi_x-symmetric}
\psi_{,x}(\rbold,\rbold_{o}; \hat{k}) &=
\dfrac{1}{\pi} \int_{\xi \in C^{+}}  \varPsi_{,x} (x,x_o;\xi)  \cos{(y \xi)} \: d\xi, 
\\
\label{eqn:IFT-phi_y-symmetric}
\psi_{,y}(\rbold,\rbold_{o}; \hat{k}) &=
-\dfrac{1}{\pi} \int_{\xi \in C^{+}} \xi  \varPsi (x,x_o;\xi) \sin{(y \xi)} \: d\xi. 
\end{align}

\subsection{Numerical approximation of the Fundamental Solution}
\label{subsec:NumericalFourierInversion}
The fundamental solution, expressed above as three indefinite IFT integrals, can not be computed analytically when water depth, and consequently wave number, change arbitrarily in one direction. In this section, the numerical aspects of its approximation and efficient numerical computation for this case are analyzed.

\subsubsection{Integration in the complex plane}

\begin{figure}
\centering
\def\svgwidth{\columnwidth}
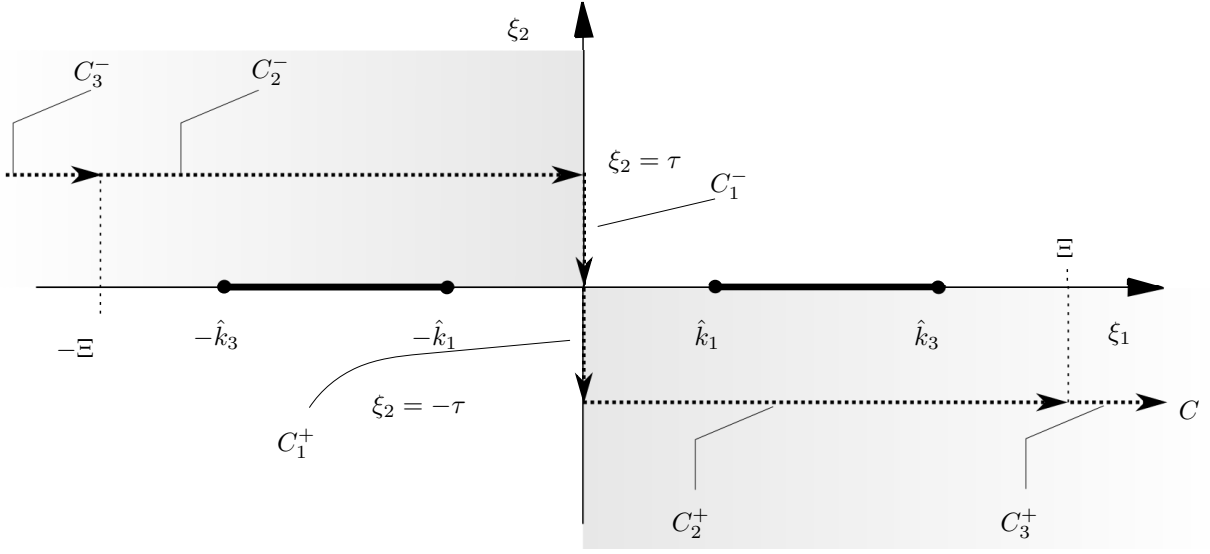
\caption{Approximation and decomposition into linear paths of the integration path in the complex $\xi$-plane used to compute the inverse Fourier transform of the function $\Psi$}
\label{fig:IntegrationPath2}
\end{figure}

For the numerical evaluation of the fundamental solution, the positive domain of integration $C^+$ is divided into three linear paths as shown in Figure \ref{fig:IntegrationPath2}, i.e., $C^{+}=C_{1}^{+} \cup C_{2}^{+} \cup C_{3}^{+}$. The first line $C_{1}^{+}= \left\lbrace \xi_{1}=0, -\tau \leq \xi_{2} \leq 0 \right\rbrace$ starts at the origin and is supposed to be very short, i.e., $\tau<<1$. The second interval $C_{2}^{+}= \left\lbrace 0 \leq \xi_{1} \leq \Xi ,\: \xi_{2}=-\tau \right\rbrace$ is finite but long enough to circumvent the roots and, finally, a third semi-infinite path $C_{3}^{+} = \left\lbrace \Xi \leq \xi_{1} ,\: \xi_{2}=-\tau \right\rbrace$ closes the domain.
Substituting this integration path in \eqref{eqn:IFT-phi-symmetric}, we compute the transformed velocity potential in the following way:
\begin{multline}
\label{eqn:IFT-phi-symmetric-NFT1}
\psi(\rbold,\rbold_o; \hat{k}) = 
\dfrac{1}{\pi} \int_{0}^{-\tau} \varPsi (x,x_o;i\xi_{2}) \cosh{(y \xi_{2})} \: d\xi _{2} +\\
\dfrac{1}{2\pi} e^{\tau y} \int_{-\Xi}^{\Xi} \varPsi (x,x_o;|\xi_{1}|-i\tau) e^{i y \xi_{1}} \: d\xi_{1} + \\
\dfrac{1}{\pi} \sinh(\tau y) \int_{0}^{\Xi} \varPsi (x,x_o;\xi_{1}-i\tau) e^{-i y \xi_{1}} \: d\xi_{1} + \\
\dfrac{1}{\pi} \int_{\Xi}^{\infty} \varPsi (x,x_o;\xi_{1}-i\tau) \cos{((\xi_{1}-i\tau)y)} \: d\xi_{1},
\end{multline}
where the first integral, along $C_{1}^{+}$, can be disregarded if $\tau$ is considered sufficiently small. The second and the third integrals correspond to the integration along the second path $C_{2}^{+}$, with the third one containing an hyperbolic sine of small argument that can also be neglected without an important loss of accuracy. The last integral corresponds to the third path $C_{3}^{+}$, where an asymptotic behavior of the integrand, defined by \eqref{eqn:FS-exponentially}, can be assumed. Under these assumptions, the approximation of the velocity potential can be finally reduced to:
\begin{multline}
\label{eqn:IFT-phi-symmetric-NFT2}
\psi(\rbold,\rbold_{o}; \hat{k}) \simeq \dfrac{1}{2\pi} e^{\tau y} \int_{-\Xi}^{\Xi} \varPsi (x,x_o;| \xi_{1} | -i\tau) e^{iy \xi_{1}} \:d\xi_{1} + \\
\dfrac{1}{2\pi} \cosh(\tau y) \Re \left\lbrace E_{1}((| x-x_o | +iy)\Xi) \right\rbrace,
\end{multline}
expression previously proposed by Belibassakis in \cite{Belibassakis2000} for the evaluation of the Green's function. However, to complete the fundamental solution kernel, we also need to compute the spatial derivatives.

The $x$-derivative of the velocity potential is obtained using the same procedure. Substituting the integration path $C^{+}$ in \eqref{eqn:IFT-phi_x-symmetric}, we have:
\begin{multline}
\label{eqn:IFT-phi_x-symmetric-NFT1}
\psi_{,x}(\rbold,\rbold_{o}; \hat{k}) = 
 \dfrac{1}{\pi} \int_{0}^{-\tau} 
 \varPsi_{,x} (x,x_o;i \xi_{2}) \cosh{(y \xi_{2})} \:d\xi_{2}
\\
+ \dfrac{1}{2\pi} e^{\tau y} \int_{-\Xi}^{\Xi} 
 \varPsi_{,x} (x,x_o;\vert \xi_{1} \vert -i\tau)  e^{i y \xi_{1}} \:d\xi_{1}
\\
+ \dfrac{1}{\pi} \sinh(\tau y) \int_{0}^{\Xi} 
 \varPsi_{,x} (x,x_o; \xi_{1} -i\tau)  e^{-i y \xi_{1}} \:d\xi_{1}
\\
+ \dfrac{1}{\pi} \int_{\Xi}^{\infty}
 \varPsi_{,x} (x,x_o; \xi_{1} -i\tau) \cos{((\xi_{1}-i\tau)y)} \:d\xi_{1},
\end{multline}
where, again, the first and the third integrals can be neglected for a small value of $\tau$. Using the asymptotic behavior of the velocity potential derivative, the integral along the semi-infinite interval $C_3^+$ can be evaluated analytically:
\begin{multline}
\int_{\Xi}^{\infty}
 \varPsi_{,x} (x,x_o; \xi_{1} -i\tau) \cos{((\xi_{1}-i\tau)y)} \:d\xi_{1}
 = \\
\dfrac{1}{2} \dfrac{e^{-| x-x_o | (\Xi-i \tau)}}{(x-x_o)^{2}+y^{2}} \left[ y \sin((\Xi-i\tau)y) -| x-x_o| \cos((\Xi-i\tau)y) \right],
\end{multline}
an substituting back in \eqref{eqn:IFT-phi_x-symmetric-NFT1} we arrive to the final approximation for the $x$-derivative of the velocity potential:
\begin{multline}
\label{eqn:IFT-phi_x-symmetric-NFT2}
\psi_{,x}(\rbold,\rbold_{o}; \hat{k})  \simeq 
\dfrac{1}{2\pi} e^{\tau y} \int_{-\Xi}^{\Xi} \varPsi_{,x} (x,x_o;| \xi_{1} | -i\tau)  e^{iy \xi_{1}} \:d\xi_{1} + \\ 
+ \dfrac{1}{2\pi} \dfrac{e^{-| x-x_o | (\Xi-i \tau)}}{(x-x_o)^{2}+y^{2}} \left[ y \sin((\Xi-i\tau)y) -| x-x_o| \cos((\Xi-i\tau)y) \right].
\end{multline}
%

However, obtaining the $y$-derivative is more involved. We see from its definition \eqref{eqn:IFT-phi_y-symmetric} that the integrand is antisymmetric, so integration in $C^{+}$ can be carried out considering only the antisymmetric part of the exponential complex function as follows:
\begin{multline}
\label{eqn:IFT-phi_y-symmetric-NFT1}
\psi_{,y}(\rbold,\rbold_{o}; \hat{k}) = 
\dfrac{1}{\pi} \int_{0}^{-\tau} \xi_{2} \varPsi (x,x_o;\xi_{2}) \sinh{(y \xi_{2})} \:d\xi _{2} +
\\
\dfrac{1}{2\pi} e^{\tau y}\int_{-\Xi}^{\Xi} 
i (\vert \xi_{1} \vert -i\tau) \varPsi (x,x_o;|\xi_{1}|-i\tau) e^{i y \xi_{1}} \:d\xi_{1} + 
\\
\dfrac{1}{\pi} \cosh{(\tau y)} \int_{0}^{\Xi} i(\xi_{1}-i\tau) \varPsi (x,x_o;\xi_{1}-i\tau) e^{-i y \xi_{1}} \: d\xi_{1} + 
\\
\dfrac{1}{\pi} \int_{\Xi}^{\infty} i(\xi_{1}-i\tau) \varPsi (x,x_o;\xi_{1}-i\tau) \sin{((\xi_{1}-i\tau)y)} \:d\xi_{1},
\end{multline}
where, once more, the first integral can be neglected, but the hyperbolic cosine that appears now in the third integral should be retained. To facilitate its numerical computation, this third integral will be extended to a symmetrical integration interval using the symmetry properties of the integrand in the following way:
\begin{multline}
\label{eqn:IFT-phi_y-symmetric-NFT2}
\int_{0}^{\Xi} i(\xi_{1}-i\tau) \varPsi (x,x_o;\xi_{1}-i\tau) e^{-i y \xi_{1}} \: d\xi_{1} =
\\
\dfrac{1}{2}\int_{-\Xi}^{\Xi} i(\vert \xi_{1} \vert -i\tau) \varPsi (x,x_o; \vert \xi_{1} \vert -i\tau) e^{-i y \xi_{1}} \:d\xi_{1} +
\\
\dfrac{1}{2}\int_{-\Xi}^{\Xi} i\xi \varPsi (x,x_o;\vert \xi_{1} \vert -i\tau) e^{-i y \xi_{1}} \:d\xi_{1},
\end{multline}
and the fourth integral, along the path $C_3^+$, is evaluated analytically using the asymptotic value of the transformed velocity potential:
\begin{multline}
\label{eqn:IFT-phi_y-symmetric-NFT3}
\int_{\Xi}^{\infty} i(\xi_{1}-i\tau) \varPsi (x,x_o;\xi_{1}-i\tau) \sin{((\xi_{1}-i\tau)y)} \:d\xi_{1} =
\\
\dfrac{1}{2} \dfrac{e^{-| x-x_o | (\Xi-i \tau)}}{(x-x_o)^{2}+y^{2}} \left[| x-x_o| \sin((\Xi-i\tau)y) +y \cos((\Xi-i\tau)y) \right],
\end{multline}
results that are substituted back in \eqref{eqn:IFT-phi_y-symmetric-NFT1} to find the final approximation for the y-derivative of the velocity potential:
\begin{multline}
\label{eqn:IFT-phi_y-symmetric-NFT4}
\psi_{,y}(\rbold,\rbold_{o}; \hat{k}) \simeq 
\dfrac{1}{2\pi} e^{\tau y}\int_{-\Xi}^{\Xi} 
i (\vert \xi_{1} \vert -i\tau) \varPsi (x,x_o;|\xi_{1}|-i\tau) e^{i y \xi_{1}} \:d\xi_{1} + 
\\
\dfrac{1}{2\pi} \cosh{(\tau y)} \int_{-\Xi}^{\Xi} i(\vert \xi_{1} \vert -i\tau) \varPsi (x,x_o; \vert \xi_{1} \vert -i\tau) e^{-i y \xi_{1}} \:d\xi_{1} +
\\
\dfrac{1}{2\pi} \cosh{(\tau y)} \int_{-\Xi}^{\Xi} i\xi \varPsi (x,x_o;\vert \xi_{1} \vert -i\tau) e^{-i y \xi_{1}} \:d\xi_{1} +
\\
\dfrac{1}{2\pi} \dfrac{e^{-| x-x_o | (\Xi-i \tau)}}{(x-x_o)^{2}+y^{2}} \left[| x-x_o| \sin((\Xi-i\tau)y) +y \cos((\Xi-i\tau)y) \right],
\end{multline}
an expression where all the integration paths are now symmetric. As we will see, this symmetry is needed for an efficient numerical evaluation using FFT.

In conclusion, equations \eqref{eqn:IFT-phi-symmetric-NFT2}, \eqref{eqn:IFT-phi_x-symmetric-NFT2} and \eqref{eqn:IFT-phi_y-symmetric-NFT4} constitute the complete kernel of the fundamental solution back-transformed to the space domain. Next step is to devise an efficient and accurate numerical procedure to evaluate these integrals.

\subsubsection{Numerical evaluation of the Fourier integrals}
The final expression of the velocity potential \eqref{eqn:IFT-phi-symmetric-NFT2} and its derivatives, \eqref{eqn:IFT-phi_x-symmetric-NFT2} and \eqref{eqn:IFT-phi_y-symmetric-NFT4}, can be calculated very efficiently by means of the FFT algorithm, as proposed in \cite{Belibassakis2000}.

Starting with the velocity potential \eqref{eqn:IFT-phi-symmetric-NFT2}, it can be expressed in a compact form as:
\begin{equation}
\psi(\rbold,\rbold_{o}; \hat{k}) \simeq \psi^{(1)} + \psi^{(2)}_{an}
\end{equation}
where $\psi^{(1)}$ is the first definite integral and $\psi^{(2)}_{an}$ represents the second analytic term on the right hand side. To compute $\psi^{(1)}$, the domain of integration is discretized using a uniform mesh of $N$ elements, where $ \xi_{l}= (l-1) \Delta \xi$ are the nodal locations of the $l=1,\cdots,N+1$ sampling points and $\Delta\xi=\Xi/N$ is the element length; in addition, the physical space $y \in [0,Y]$ is discretized with a similar uniform distribution of nodes $y_{j}=(j-1)\Delta y$, for $j=1,\cdots,N+1$, separated a distance $\Delta y=Y/N$ using a fixed value of $\Delta y = \pi/\Xi$. The symmetry property expressed in \eqref{eqn:symetry-xi} assures that the discrete values of the transformed velocity potential $\varPsi_{l}=\varPsi(x,x_o; \xi_{l}-i\tau)$, for $l=1,\ldots,N+1$, are symmetric with respect to $\xi_{1}$ and hence $\varPsi_{l}=\varPsi_{2N-l+2}$, for $l=2,\ldots,N$. Based on these discretizations and applying the IFFT algorithm, the value of $\psi^{(1)}$ can then be approximated by the finite series:
\begin{equation}
\psi^{(1)}=\dfrac{1}{2\pi} e^{\tau y_{j}} \left[ \sum_{l=1}^{M} \varPsi_{l} e^{i\frac{2\pi}{M}(j-1)(l-1)}
\right]
\end{equation}
where the total number of points $M=2N$ is selected as a power of two to be efficiently evaluated.

In a similar way, the $x$-derivative given by \eqref{eqn:IFT-phi_x-symmetric-NFT2} can be written as the addition of two terms:
\begin{equation}
\psi_{,x}(\rbold,\rbold_{o}; \hat{k}) \simeq \psi^{(1)}_{,x} + \psi^{(2)}_{,x \: an}
\end{equation}
with a definite integral $\psi^{(1)}_{,x}$ that needs to be evaluated at the same interval. Noting that the discrete values of the velocity potential derivative $\varPsi_{,x \: l}=\varPsi_{,x}(x,x_o; \xi_{l}-i\tau)$, for $l=1,\ldots,N+1$, are symmetric with respect to $\xi_{1}$ and using the same discretization, we can sample $\varPsi_{,x \: l}=\varPsi_{,x}(x,x_o; \xi_{l}-i\tau)$, for $l=1,\ldots,N+1$, and apply the symmetry property $\varPsi_{,x \: l}=\varPsi_{,x \: 2N-l+2}$, for $ l=2,\ldots,N$, to approximate $\varPsi_{,x}^{(1)}$ as the IFFT sequence:
\begin{equation}
\psi_{,x}^{(1)}=\dfrac{1}{2\pi} e^{\tau y_{j}} \left[ \sum_{l=1}^{M} \varPsi_{,x \: l} e^{i\frac{2\pi}{M}(j-1)(l-1)}
\right].
\end{equation}

Finally, according to equation \eqref{eqn:IFT-phi_y-symmetric-NFT4}, the $y$-derivative of the velocity potential can be decomposed into four different terms:
\begin{equation}
\psi_{,y}(\rbold,\rbold_{o}; \hat{k}) \simeq \psi^{(1)}_{,y} + \psi^{(2)}_{,y} + \psi^{(3)}_{,y} + \psi^{(4)}_{,y \: an}
\end{equation}
with three definite integrals and one analytical term. The complex parameter present in the integrands of $\psi^{(1)}_{,y}$ and $\psi^{(2)}_{,y}$ is $\xi^{s} = \vert \xi_{1} \vert -i\tau$, a symmetric function with respect to variable $\xi_{1}$. This variable can be defined in a discrete form as $\xi^{s}_{l}= (l-1) \Delta \xi -i \tau$, for 
$l=1,\ldots,N+1$, and $\xi^{s}_{2N-l+2}=\xi^{s}_{l}$, for  $l=2,\ldots,N$. By doing that, the first integral $\psi^{(1)}_{,y}$, can be evaluated applying the IFFT algorithm as we did before:
\begin{equation}
\psi_{,y}^{(1)}=\dfrac{1}{2\pi} e^{\tau y_{j}} \left[ \sum_{l=1}^{M} i\xi^{s}_{l} \varPsi_{l} e^{i\frac{2\pi}{M}(j-1)(l-1)}
\right]
\end{equation}
and, on the contrary, the second integral $\psi^{(2)}_{,y}$ is approximated using the FFT algorithm:
\begin{equation}
\psi_{,y}^{(2)}=\dfrac{1}{2\pi} \cosh{(\tau y_{j})} \left[ \sum_{l=1}^{M} i\xi^{s}_{l} \varPsi_{l} e^{-i\frac{2\pi}{M}(j-1)(l-1)}
\right]
\end{equation}
due to its negative exponential term.

In the last integral $\psi^{(3)}_{,y}$, the complex parameter $\xi$ is defined in discrete form as
$\xi^{a}_{l}= (l-1) \Delta \xi -i \tau$, for 
$l=1,\ldots,N+1$,
and 
$\xi^{a}_{2N-l+2}=-\xi^{a}_{l}$, for  $l=2,\ldots,N$.
Substituting and performing the FFT, the last definite integral is approximated as follows:
\begin{equation}
\psi_{,y}^{(3)}=\dfrac{1}{2\pi} \cosh{(\tau y_{j})} \left[ \sum_{l=1}^{M} i\xi_{l}^{a} \varPsi_{l} e^{-i\frac{2\pi}{M}(j-1)(l-1)}
\right]
\end{equation}
closing the derivation of a complete fundamental solution kernel.

The accuracy of this approach highly depends on a proper selection of parameters $\tau$, $\Xi$ and $N$. As explained by Belibassakis \cite{Belibassakis2000}, for the calculation of $\psi(\rbold,\rbold_{o}; \hat{k})$ using this technique, the value of $\tau$ must be small enough to make it possible to neglect the contribution of the first and third integrals of \eqref{eqn:IFT-phi-symmetric-NFT1} and \eqref{eqn:IFT-phi_x-symmetric-NFT1}, and the first term of \eqref{eqn:IFT-phi_y-symmetric-NFT1}, but at the same time, it can not be too small because the aliasing effect is attenuated at least by factor of $\exp(-2\tau Y)$. On the other hand, experience demonstrates that a value of $\Xi \approx 4-6\hat k^{*}$, being $\hat k^{*}$ the maximum value of $\hat k(x)$, is large enough to make the asymptotic expression \eqref{eqn:FS-exponentially} valid, and consequently the approximation of integrals along the path $C_{3}^{+}$. In our calculations, we have used a sampling of $M=2^{12}=4096$ points inside the interval $[-\Xi,\Xi]$, fixing the other two parameters to $\Xi=6\hat k^{*}$ and $\tau=\Delta \xi$.

\subsubsection{FEM solution of the transformed velocity potential}
As we have seen, in order to evaluate the fundamental solution, it is necessary to solve the transformed one-dimensional wave transmission problem defined by equations (\ref{eqn:1Dvarphi in 2}-\ref{eqn:bc-varphi-x}) for different values of the Fourier parameter $\xi$. For this task, Belibassakis \cite{Belibassakis2000} proposes a second-order central finite difference scheme. In our experience, the use of the finite element method improves the solution near the source point $x_{o}$, increasing this way the final accuracy of the fundamental solution.

Applying the method of weighted residuals, with a test function $w(x)$ defined in the domain $[a,b]$, the weak form of equation \eqref{eqn:1Dvarphi in 2} can be expressed:
\begin{equation}
\label{eqn:FEM1Dweakform1}
\int_{a}^{b} w( \varPsi_{,xx} + \kappa^{2}(x)\varPsi +  \delta(x-x_o)) \: dx = 0,
\end{equation}
and integrating by parts:
\begin{equation}
\label{eqn:FEM1Dweakform3}
\int_{a}^{b}[w_{,x} \varPsi_{,x}-w\kappa^{2}(x)\varPsi] \: dx = 
w(x_o) + i\beta(\xi) w(b) \varPsi(b) + i\alpha(\xi) w(a) \varPsi(a),
\end{equation}
where we have substituted the matching conditions \eqref{eqn:bc-varphi-x} that close the transformed domain.

Using the classical finite element Galerkin formulation, the transformed velocity potential and the weighted residual function are approximated as:
\begin{equation}
\label{eqn:FEM-approach}
\varPsi(x)= \sum_{j=1}^{n} N_{j}(x)\varPsi_{j}, \quad w(x) = \sum_{j=1}^{n} N_{j}(x) w_{j}
\end{equation}
where  $\varPsi_{j}$ and $w_{j}$ are the corresponding nodal values, $N_j(x)$ are linear shape functions and $n$ is the number of nodes distributed in the domain. Substituting the discretization \eqref{eqn:FEM-approach} into \eqref{eqn:FEM1Dweakform3}, the following FEM system is obtained:
\begin{equation}
\label{eqn:FEM-SYSTEM}
\left[
\begin{matrix}
a_{11}-i\alpha(\xi) & a_{12} & \cdots  &  \cdots  & 0 \\
a_{21} & a_{22} & \cdots  &  \cdots  & 0 \\
\vdots & \vdots & a_{ij}  &  \ddots  & \vdots\\
\vdots & \vdots &  \ddots &  \ddots  & \vdots\\
0 & 0 & \cdots  &  \cdots  & a_{n n}-i\beta(\xi)
 \end{matrix}
 \right]
 \left\lbrace 
 \begin{matrix}
 \varPsi_1 \\
 \vdots \\
 \varPsi_j\\
 \vdots \\
 \varPsi_{n}
 \end{matrix}
 \right\rbrace
 =
 \left\lbrace
 \begin{matrix}
 0\\
 \vdots\\
 1 \\
 \vdots\\
 0
 \end{matrix}
 \right\rbrace,
 \end{equation}
where the components of the tridiagonal matrix:
\begin{equation}
a_{ij} = \int_{a}^{b} \left[  N_{i,x}N_{j,x} - \kappa^{2}(x)  N_i N_j\right] \: dx.
\end{equation}
are computed with a two-point Gauss quadrature.

The right hand side of system \eqref{eqn:FEM-SYSTEM} contains the contribution of the load, where a node $j$ is purposely located at $x_o$. As a practical rule, we use at least 20 elements per wave-length to discretize the transformed one-dimensional transmission problems.

\subsection{Validation of the fundamental solution for constant water depth}
To check the accuracy of the numerical integration process described in Section \ref{subsec:NumericalFourierInversion} to approximate the MSE fundamental solution, we solve the problem for constant water-depth and compare the numerical results with the analytical solution of the equivalent Helmholtz problem for constant wave number given by:
\begin{equation}
\psi_{H}=  \dfrac{i}{4} H_{0}^{(1)}(kr) , \quad
\psi_{H,r} = - \dfrac{i}{4} k H_{1}^{(1)}(kr),
\end{equation}
where $H_{0}^{(1)}$ and $H_{1}^{(1)}$ are Hankel functions of the first kind of order zero and one, $r=\vert \rbold-\rbold_o\vert$ is the distance in the radial direction and subindex $H$ refers to Helmholtz solution. Using the series expansion of the Bessel's function for small argument \cite{Abramowitz} the fundamental solution $\psi_H$ near the source can be expressed as:
\begin{equation}
\psi_{H}^p = -\frac{1}{2\pi} (\ln\frac{kr}{2}+\gamma) + \frac{i}{4} + \mathcal{O} \left( (kr)^2\ln(kr) \right)
\label{eqn:seriesPhiH}
\end{equation}
where $\gamma$ is the Euler-Mascheroni constant and the superindex $p$ denotes the polynomial approximation. The expansion contains a weak singularity of the the real part and a constant imaginary value.

\begin{figure}
\begin{center}
\includegraphics[width=0.8\columnwidth]{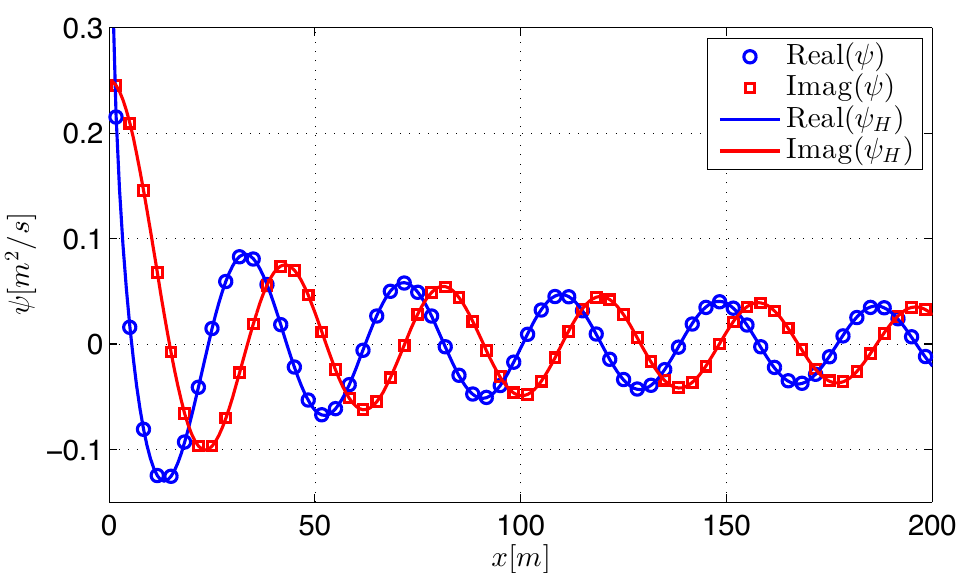} \\
\includegraphics[width=0.8\columnwidth]{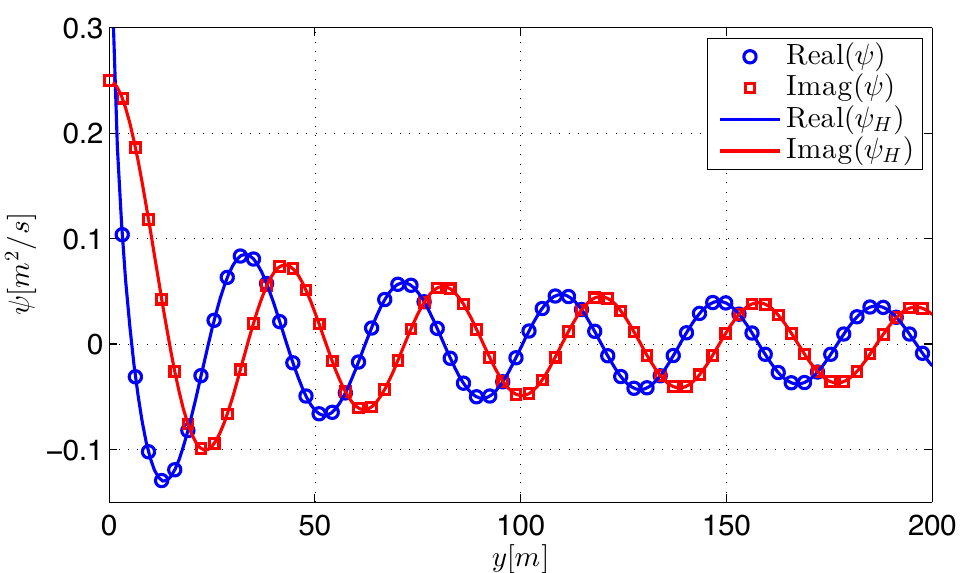}
\caption{Comparison of analytical and numerical fundamental solution $\psi$ for constant water depth. Cross section of the 2D fundamental solution along the lines $y=0$ (top) and $x=0$ (bottom)}
\label{fig-FS_XY}
\end{center}
\end{figure}

\begin{figure}[t]
\begin{center}
\includegraphics[width=0.8\columnwidth]{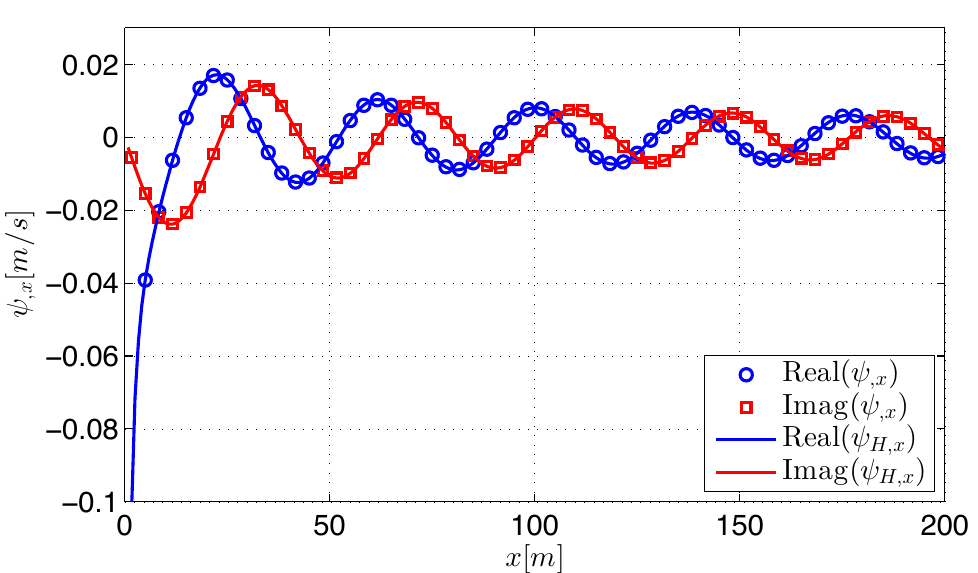} \\
\includegraphics[width=0.8\columnwidth]{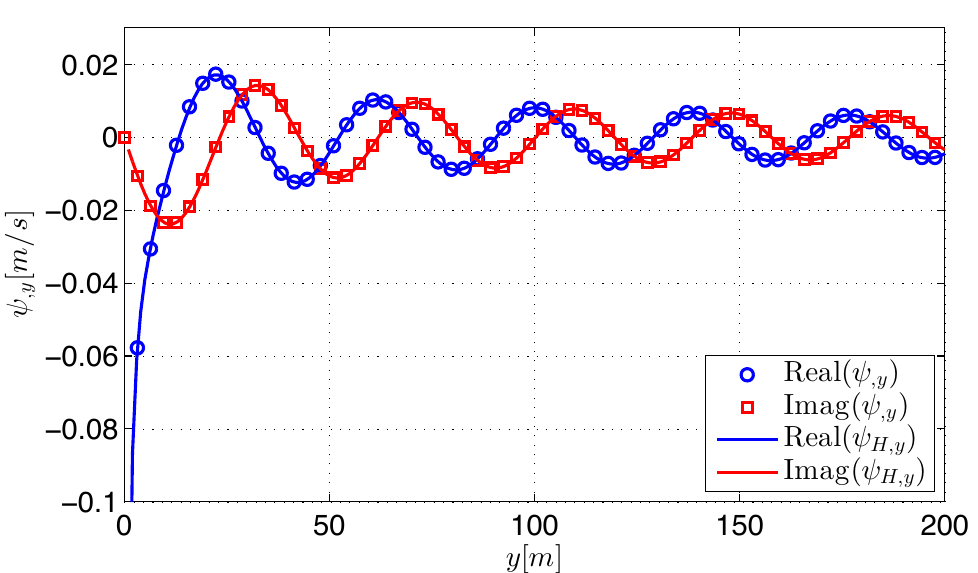}
\caption{Comparison of analytical and numerical fundamental solution derivatives $\psi_{,x}$ and $\psi_{,y}$ for constant water depth. Cross section of the $x$-derivative along the line $y=0$ (top) and the $y$-derivative along $x=0$ (bottom)}
\label{fig-DFS_XY}
\end{center}
\end{figure}

Figures \ref{fig-FS_XY} and \ref{fig-DFS_XY} show a comparison between the numerical fundamental solution $\psi$ and the analytical solution $\psi_{H}$ for constant water depth. The wave-period and water-depth used are $T=5s$ and $h=14m$ respectively. The solution profiles at $y=0$ and $x=0$ shown in Figure \ref{fig-FS_XY} are in excellent agreement. The same degree of approximation is obtained for $x$ and $y$ derivatives, represented in Figure \ref{fig-DFS_XY}. 

\section{Boundary element formulation}
\label{sec-BEM}
In this section we present the basis of the BEM for the Helmholtz problem, originated from an integral equation expressing a reciprocity relation between the unknown velocity potential field $(\hat{\phi},\nabla \hat{\phi})$ and the fundamental solution kernel $(\psi,\nabla \psi)$.

\subsection{Direct boundary integral equation}
The BEM formulation for an acoustic medium is well known and can be found in different texts \cite{Wu2000,Atalla2015}. Let us consider a domain $\Omega$ closed by a boundary $\Gamma$ of outward normal $\nbold$. Multiplying both sides of Helmholtz equation \eqref{eqn:Helmholtz} by the fundamental solution and applying Green's second identity and Sommerfeld's radiation condition, the following Boundary Integral Equation (BIE) is obtained in absence of internal loads:
\begin{equation}
C(\rbold_o) \hat{\phi}(\rbold_o)  
+ 
\int_{\Gamma}  \nabla \psi(\rbold,\rbold_o;\hat{k}) \cdot \nbold \: \hat{\phi}(\rbold) \:d\Gamma 
-
\int_{\Gamma} \psi(\rbold,\rbold_o;\hat{k}) \nabla \hat{\phi}(\rbold) \cdot \nbold \:d\Gamma 
= 0,
\label{eqn:BEMBoundaryIntegralEquation}
\end{equation}
where $\nabla \hat{\phi} \cdot \nbold = {\partial \hat{\phi}}/{\partial n} =\hat{q}$ is the normal flux, $\rbold_o$ is the collocation point and $C(\rbold_o)$ is a geometrical coefficient, function of the regularity of the boundary, that takes the value $C(\rbold_o) ={\theta(\rbold_o)}/{2\pi}$ with $\theta(\rbold_o)$ being the internal angle of the boundary at point $\rbold_o$. 

This BIE can be solved numerically discretizing the boundary into $n_e$ elements, $\Gamma = \bigcup\limits_{i=1}^{n_e} \Gamma_{e}$ and approximating the fields of modified velocity potential $\hat{\phi}(\rbold)$ and normal flux $\hat{q}(\rbold)$ using isoparametric linear elements:
\begin{equation}
\label{eqn:BEMdiscretizacion}
\hat{\phi}(\zeta)=\sum_{i=1}^{n} N_{i}(\zeta) \hat{\phi}_{i} ,
\hspace{0.5cm}
\hat{q}(\zeta)=\sum_{i=1}^{n} N_{i}(\zeta) \hat{q}_{i},
\end{equation}
where $\hat{\phi}_{i}$ and $\hat{q}_{i}$ represent the nodal values of velocity potential and flux, $N_{j}(\zeta)$ are the element shape functions and $n$ is the number of nodes per element. Introducing this approximation in \eqref{eqn:BEMBoundaryIntegralEquation}, the discretized form of the BIE can be written as:
\begin{equation}
\label{eqn:BEMdiscretized BIE}
C_{i} \hat{\phi}_{i} + \sum_{e=1}^{n_{e}} \int_{\Gamma_{e}} \nabla \psi(\rbold,\rbold_i;\hat{k}) \cdot \nbold N_{j} \: \hat{\phi}_{j} \:d\Gamma_{e} = 
\sum_{e=1}^{n_{e}} \int_{\Gamma_{e}}  \psi(\rbold,\rbold_i;\hat{k}) N_{j} \hat{q}_{j} \:d\Gamma_{e},
\end{equation}
an equation that can be expressed in matrix form as:
\begin{equation}
\label{eqn:hatHphatGq}
\hatHbold \hatphibold = \hatGbold \hatqbold,
\end{equation}
being $\hatphibold$ the nodal vector of modified velocity potentials $\hat{\phi}_j$, $\hatqbold$ the nodal vector of fluxes $\hat{q}_j$, together with the matrix coefficients:
\begin{align}
\hat{H}_{ij} & = C_{i}\delta_{ij} + \sum_{e=1}^{n_{e}} \int_{\Gamma_{e}} \nabla \psi(\rbold,\rbold_i;\hat{k}) \cdot \nbold \: N_{j} \:d\Gamma_{e} \label{eqn:Hij} \\
\hat{G}_{ij} & = \sum_{e=1}^{n_{e}} \int_{\Gamma_{e}}  \psi(\rbold,\rbold_i;\hat{k}) N_{j} \:d\Gamma_{e} \label{eqn:Gij}
\end{align}
where $\delta_{ij}$ is the Kronecker $\delta$-function and $\rbold_i$ the position of node $i$.

\begin{figure}
\centering
\includegraphics[width=0.8\columnwidth]{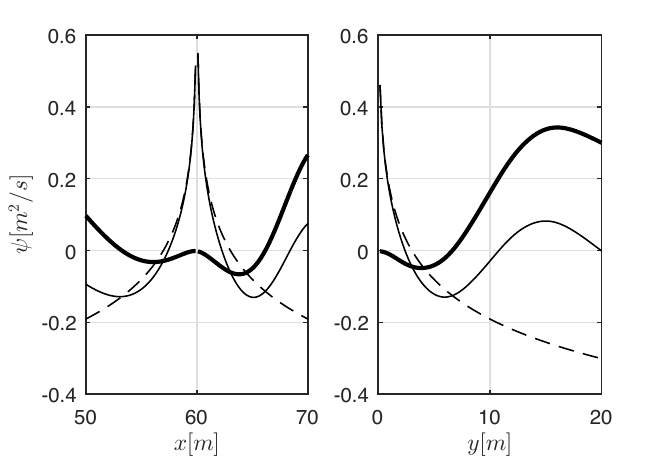}
\caption{Behavior of the Green's function $\psi(\rbold,\rbold_o;\hat{k})$ (thin solid line), the kernel $\wp(\hat{k}_o r)$ (dashed line) and the regularized integrand (thick solid line). The source point is located at $x_o=60m$ and the bathymetry is defined in Section \ref{SUBSECTION: NM CANAL} with a wave period $T=5s$}
\label{fig: el phase field}
\end{figure}

In general, the boundary element integrals present in \eqref{eqn:Hij} and \eqref{eqn:Gij} can be computed using a standard Gauss quadrature formula. But when the collocation point is located in one of the element nodes, integral \eqref{eqn:Gij} becomes weakly singular and needs a special treatment. In that particular case, it is possible to express the integrand as the sum of two terms:
\begin{equation}
\hat{G}_{ij} = \int_{\Gamma_{e}}  [ \psi(\rbold,\rbold_i;\hat{k}) - \wp(\hat{k}_i r) ] N_{j} \:d\Gamma_{e} + \int_{\Gamma_{e}} \wp(\hat{k}_i r) N_{j} \:d\Gamma_{e}
\label{eqn:GijRegular}
\end{equation}
one completely regular treated by a standard Gauss quadrature and another term of order $\mathcal{O}(\ln r)$ that can be evaluated numerically using a special quadrature. We have used the real part of the first term of the series expansion of $\psi_H$ around $\rbold_i$, obtained in \eqref{eqn:seriesPhiH}, to define the kernel:
\begin{equation}
\wp(\hat{k}_i r) = -\frac{1}{2\pi} (\ln\frac{\hat{k}_i r}{2}+\gamma)
\label{eqn:fRegular}
\end{equation}
needed to regularize the first integral, see Figure \ref{fig: el phase field}. Details of the numerical treatment of the second integral can be found in \cite{Dominguez1993}.

Finally, back-substituting the change of variable introduced in \eqref{eqn:ChangeVariable}, the BEM system \eqref{eqn:hatHphatGq} can be rewritten in terms of the original velocity potential as:
\begin{equation}
\label{eqn:HpGq}
\Hbold \phibold = \Gbold \bold{q},
\end{equation}
a complex non-symmetrical linear system that provides the solution of the problem.

\subsection{Integral formulation for scattering problems}
In wave transmission problems, it is usually interesting to consider the effect on the object under study of an incident wave radiated from a distant source. In these scattering problems \cite{Wu2000}, we divide the total velocity potential into a scattered wave and an incident wave:
\begin{equation}
\hat\phi(\rbold) = \hat\phi_{sc}(\rbold) + \hat\phi_{in}(\rbold)
\label{eqn:DecompositionScattering}
\end{equation}
where $\hat\phi_{sc}$ is the scattered field radiated by the object that should satisfy the Helmholtz equation and $\hat\phi_{in}$ is the incident field that would exist in the absence of obstacles. If we substitute this decomposition in \eqref{eqn:BEMBoundaryIntegralEquation}, the BIE adopts the new form:
\begin{multline}
C(\rbold_o) \hat{\phi}(\rbold_o)  
+
\int_{\Gamma}  \nabla \psi(\rbold,\rbold_o;\hat{k}) \cdot \nbold \: \hat{\phi}(\rbold) \:d\Gamma
-  
\int_{\Gamma} \psi(\rbold,\rbold_o;\hat{k}) \: \nabla \hat{\phi}(\rbold) \cdot \nbold \:d\Gamma
 = 
\hat\phi_{in}(\rbold_o).
\label{eqn:BEMBoundaryIntegralEquationScattering}
\end{multline}

Using the same approximation described in \eqref{eqn:BEMdiscretizacion} for the modified velocity potential and the normal fluxes on the boundary, the discrete matrix form of the BEM for scattering problems becomes:
\begin{equation}
\Hbold \phibold - \Gbold \bold{q} = \phibold_{in}
\label{eqn:BEMBoundaryIntegralEquationScatteringMatrix}
\end{equation}
where the change of variable \eqref{eqn:ChangeVariable} has been applied and $\phibold_{in}$ is a free-term vector containing the evaluation of the incident potential $\phi_{in}(\rbold)$ at the nodal positions. Therefore, the only difference with the original BEM system \eqref{eqn:HpGq} is in the free term due to the incident wave.
\section{Numerical examples}
\label{SECTION:NUMERICAL EXAMPLES}
In this section, different numerical examples are presented in order to demonstrate the potential of the proposed BEM formulation for the solution of water-wave transmission problems in variable bathymetries. The first and the second problem present analytical solution and are used to test the accuracy of the fundamental solution and the BEM scheme. The variable used to compare with the analytical solution is the wave amplification factor WAF, a common design parameter in engineering applications that is defined as the ratio between the local value of the wave height and a reference value, normally selected as the incident wave height in deep water. Due to the linear relation between the water surface elevation and the velocity potential, i.e, $|\phi|=gH/(2\omega)$, the WAF is also the ratio between the absolute value of the velocity potential and the absolute value of the incident velocity potential in deep water.

The velocity potential for an incident wave traveling in an infinite two-dimensional domain with unidirectional variable bathymetry, under the mild-slope assumption in conjunction with very slowly varying bathymetry and discarding strong diffraction effects, is given by the expression \cite{PanchangChen2000}:
\begin{equation}
\label{eqn:NE Incident Wave}
\phi_{in}(\rbold)= \vert \phi_{o} \vert \, A(x) \, \exp{ \left[ ik(x)y \sin{\theta + \int_{x_0}^{x} ik(\eta) \cos\theta \: d\eta} \right]}
\end{equation} 
where wave amplitude $A(x)$ is a function containing the shoaling and refraction coefficients and $\theta$ is the angle between the incident wave and the direction of variable wave number. This solution is a necessary ingredient of the BEM formulation for scattering problems \eqref{eqn:BEMBoundaryIntegralEquationScatteringMatrix}.

\subsection{Shoaling effect in a channel}
\label{SUBSECTION: NM CANAL}
The shoaling effect can be observed under stationary conditions in waves traveling from deep to shallow waters. When the waves arrive to the shallow water they slow down, the wave length is gradually reduced and, because the energy flux must remain constant, a reduction in the group velocity is compensated by an increase in the wave height.

\begin{figure}
\centering
\def\svgwidth{0.9\columnwidth}
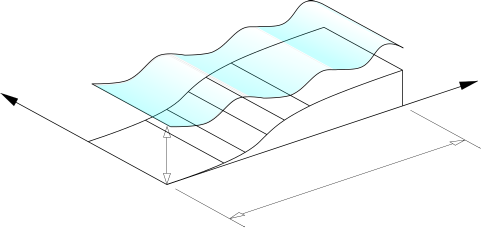 \\
\text{(a)} \\
\def\svgwidth{0.9\columnwidth}
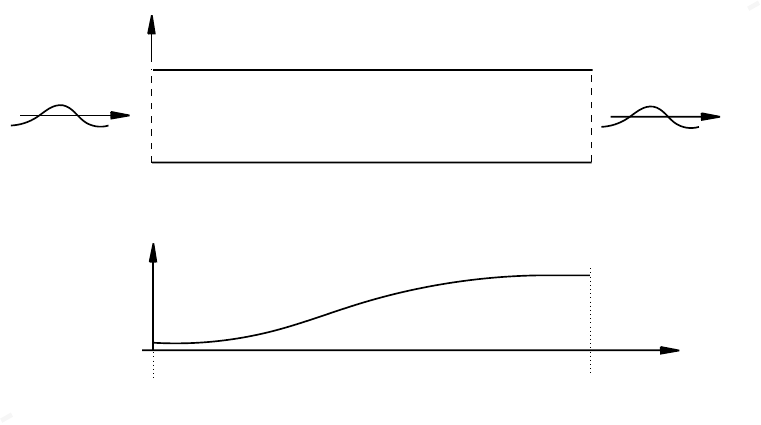 \\
\text{(b)}
\caption{Channel with variable water depth $h(x)$ in the longitudinal direction. Model of the closed rectangular domain (a). Boundary conditions and associated wave-number $k(x)$ in the x-direction (b)}
\label{fig:CANAL_SCHEME}
\end{figure}

To reproduce this phenomena, we model a rectangular channel of length $L=70m$ in the x-direction that is discretized with a uniform mesh of linear boundary elements. The lateral walls present zero normal-flux conditions, the incident velocity potential $\phi_{in}$ is imposed as boundary condition at $x=0m$ and $x=70m$. The period of the incoming water wave is $T=5s$ and its wave-length at the entry point is $L_{0}=39m$. A sketch of the configuration is represented in Figure \ref{fig:CANAL_SCHEME}.

Water depth $h(x)$ is supposed to decrease monotonically along the channel from $14m$ to $0.5m$ between $x=0m$ and $x=70m$. This means that, for the considered initial wave-length of 39 m, waves travel from intermediate water-depths to shallow waters. In the transition zone, the water depth function is mathematically approximated by a cubic polynomial:
\begin{equation}
\label{eqn:NE:Canal:profile}
h(x)=a_0+a_1x+a_2x^{2}+a_3x^{3}
\end{equation}
with $a_0=14m$, $a_1=0m^{-1}$, $a_2=-8.2653\times 10^{-3}m^{-2}$ and $a_3=7.8717\times 10^{-5}m^{-3}$. The associated wave number $k(x)$ for this water-depth function is evaluated using the dispersion relation \eqref{eqn:dispersionEquation} for the fixed frequency of the incident wave.

\begin{figure}
\begin{center}
\includegraphics[width=0.8\columnwidth]{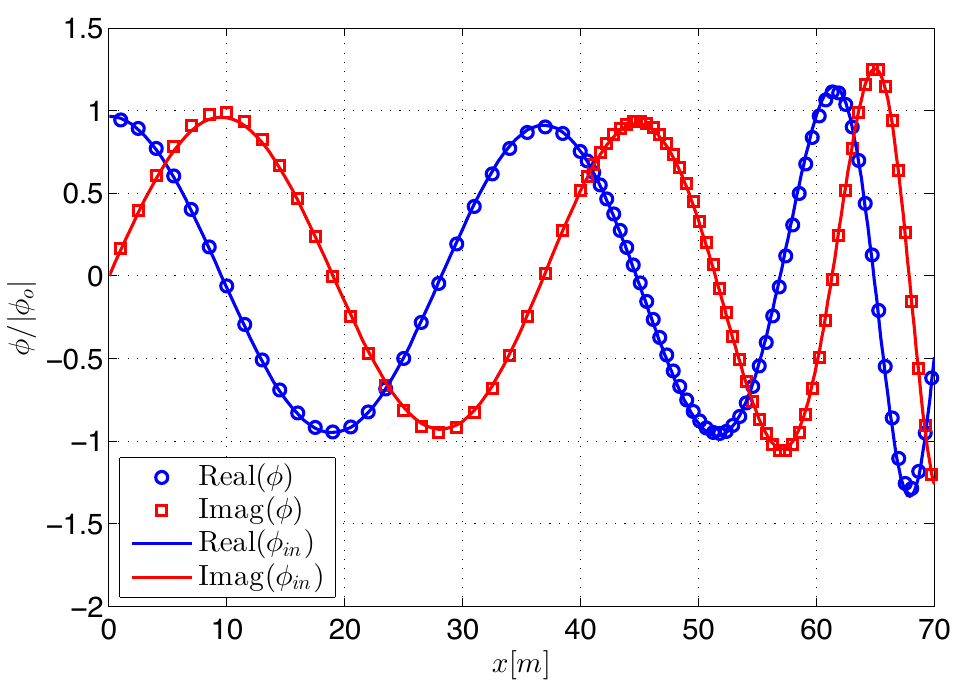}
\caption{Shoaling effect in a channel with monotonically decreasing water depth. Real and imaginary parts of the wave amplification factor obtained using BEM (dots) compared with the analytical solution (solid-line) given by equation \eqref{eqn:NE Incident Wave}}
\label{fig-Canal}
\end{center}
\end{figure}

The variation of the velocity potential along the channel is obtained at the lateral walls. Figure \ref{fig-Canal} shows a comparison of the WAF obtained with BEM and the analytical solution for the incident wave given by equation \eqref{eqn:NE Incident Wave}. Both solutions present very good agreement, demonstrating the ability of the proposed boundary element formulation to represent this phenomenon and to compute accurate shoaling coefficients in variable water depths.

\subsection{Scattering by a cylinder with variable bathymetry}
\label{SUBSECTION: NE CYLINDER CONSTANT BATHYMETY}
\begin{figure}
\centering
\def\svgwidth{0.8\columnwidth}
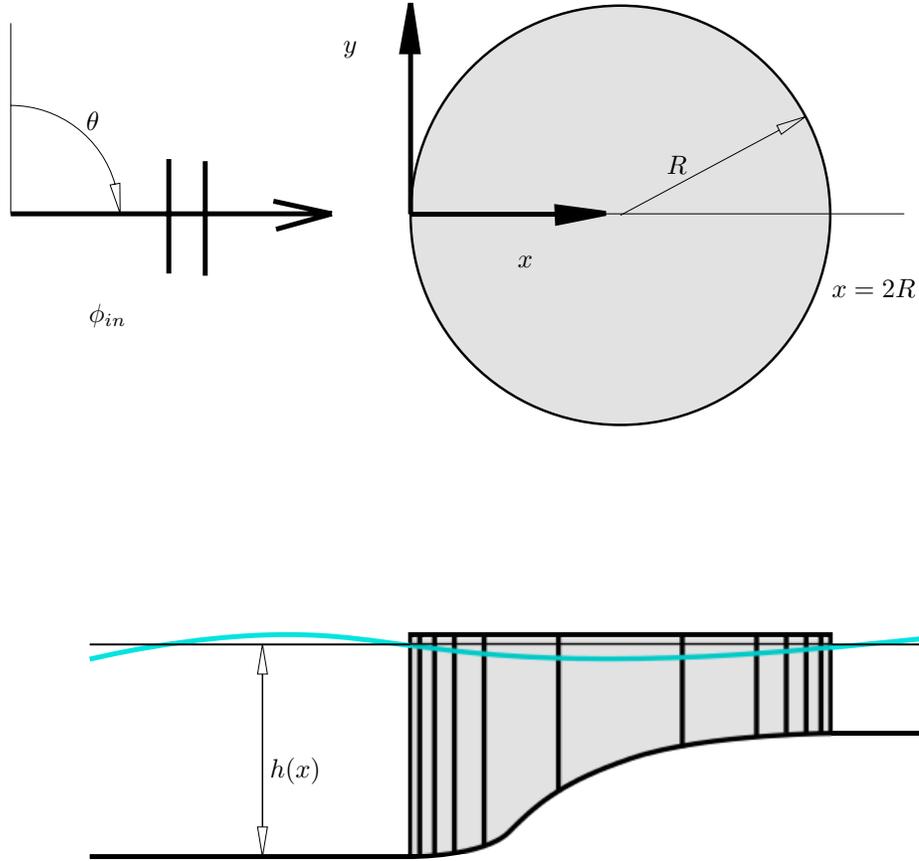
\caption{Rigid circular cylinder mounted on waters with variable depth $h(x)$. Top and side views of the cylinder with an incident wave $\phi_{in}$ of incidence angle $\theta$ relative to the bathymetric lines}
\label{fig:CYLINDER_SCHEME}
\end{figure}

In this section we study the wave scattering produced by a cylinder in waters of variable depth. A rigid cylinder of radius $R=25 m$ is fixed on a seabed of depth $h(x)$ varying in one direction, as represented in Figure \ref{fig:CYLINDER_SCHEME}. Plane waves of potential $\phi_{o}$ are incident from infinity with wave period $T$ and incidence angle $\theta$. To model the open domain, the boundary of the cylinder is uniformly discretized using 320 linear BEM elements of the same size.

\begin{figure}
\begin{center}
\subfloat[]{\includegraphics[width=7cm]{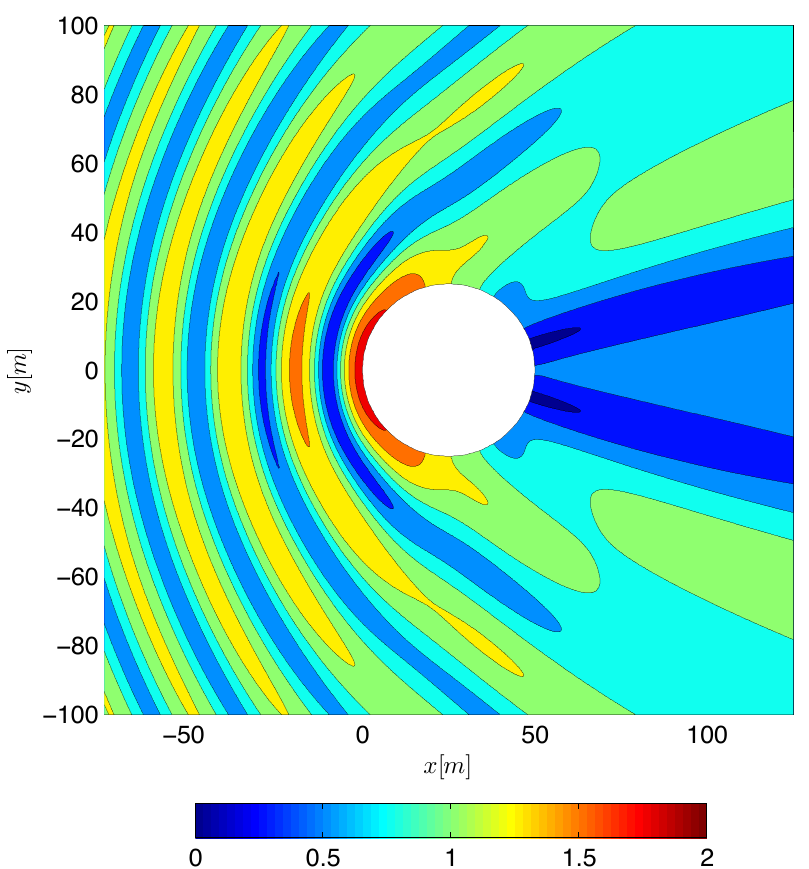}}
\vspace{0.05cm}
\subfloat[]{\includegraphics[width=7cm]{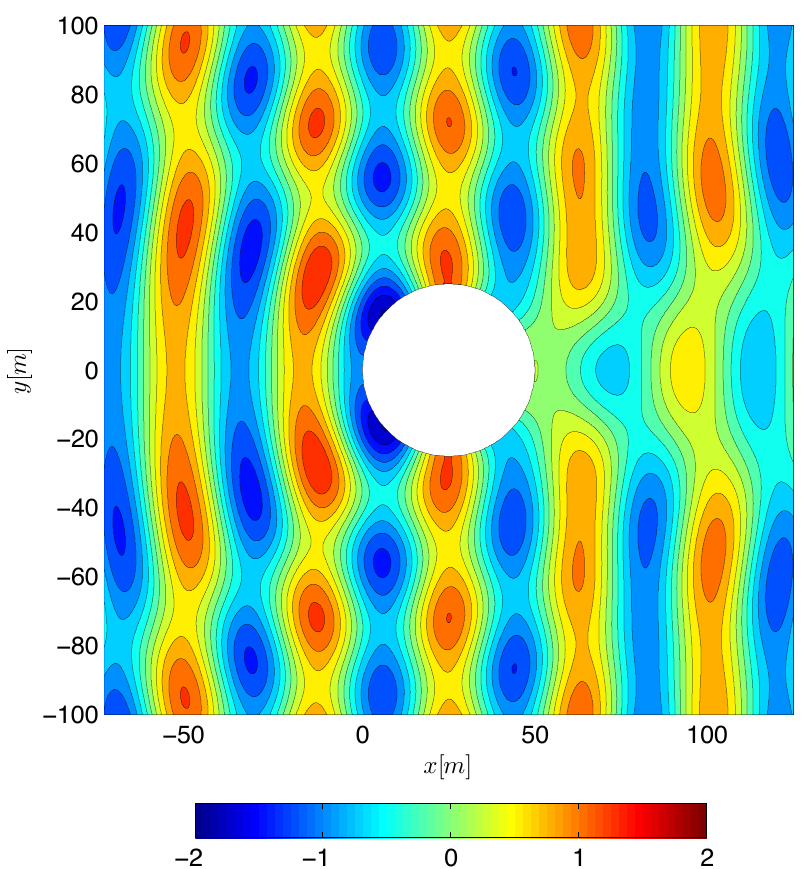}}
\caption{Interference patterns obtained using BEM for a cylinder in constant water depth $h(x)=14 m$ and an incident wave of period $T=5s$. Contours of magnitude (a) and real part (b) of the normalized velocity potential $\phi / \vert\phi_{o}\vert$}
\label{fig:Cylinder Constant Bathymetry}
\end{center}
\end{figure}

\begin{figure}
\begin{center}
\includegraphics[width=0.8\columnwidth]{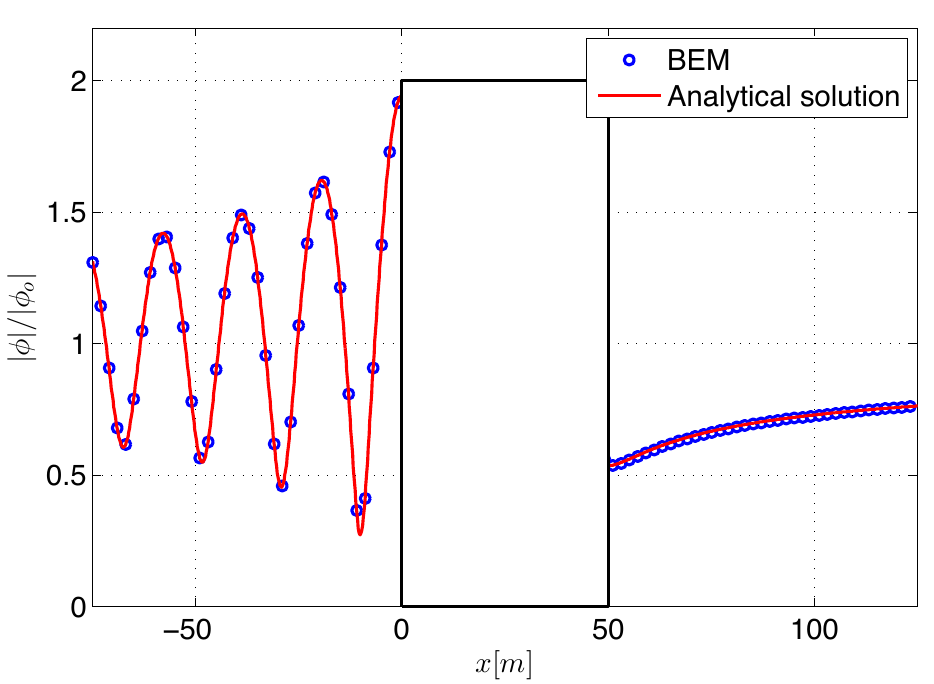}
\caption{Wave scattering by a cylinder in constant water depth. Comparison with the analytical solution of the normalized magnitude of the velocity potential $\vert\phi\vert / \vert\phi_{o}\vert$ at $y=0$ obtained using BEM}
\label{fig:Cylinder Constant Bathymetry PROFILE}
\end{center}
\end{figure}

We will apply first the proposed BEM scheme \eqref{eqn:BEMBoundaryIntegralEquationScatteringMatrix} to the particular case of constant bathymetry in order to check the accuracy achieved by the numerical fundamental solution. The wave scattering by a circular cylinder in an infinite homogeneous medium is a well known problem with analytical solution obtained by McCamy and Fuchs \cite{MacCamy1954} that is commonly used to validate numerical algorithms for diffraction problems \cite{Berkhoff1976,Mei1983}. Considering a constant water depth $h(x)=14 m$, the dispersion relation \eqref{eqn:dispersionEquation} yields for this case a wave number $k=0.42 m^{-1}$. The magnitude and real part of the normalized velocity potential $\phi(x,y)/\vert\phi_{o}\vert$ computed with BEM are represented in the form of contour plots in Figure \ref{fig:Cylinder Constant Bathymetry}. As expected, we obtain a diffraction pattern that decays inversely with the square root of the radial distance and observe a strong shadow region in the rear part of the cylinder. Accuracy of this solution is demonstrated in Figure \ref{fig:Cylinder Constant Bathymetry PROFILE}, where the normalized velocity potential is compared with the analytical solution for a cross section located at $y=0$.

\begin{figure}
\begin{center}
\includegraphics[width=8cm]{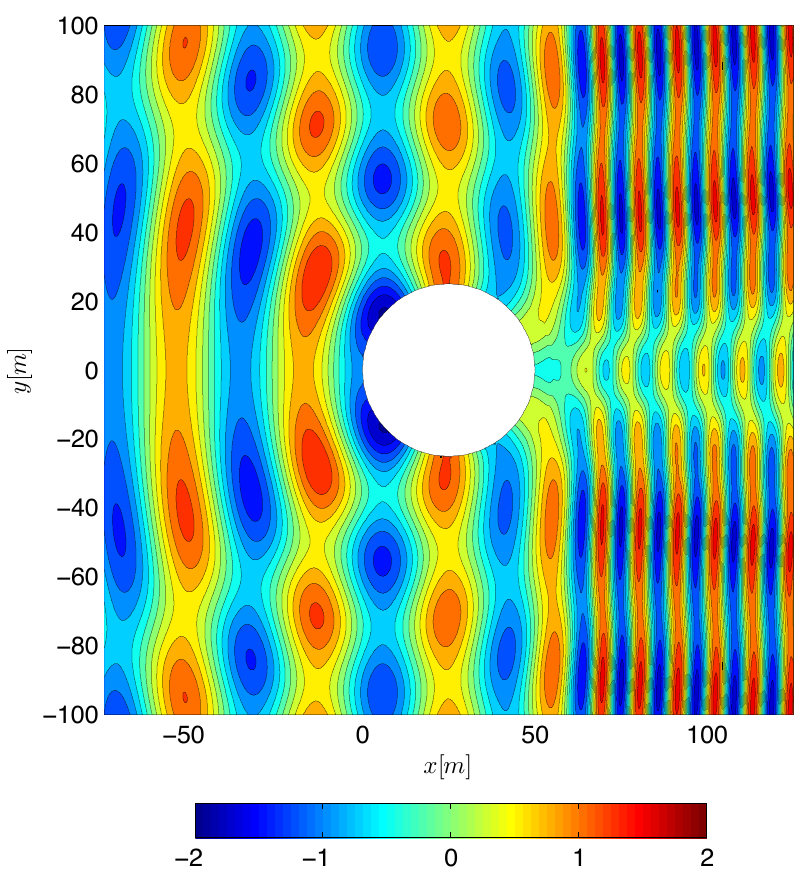}
\caption{BEM solution for a decreasing water depth in the $x$-direction and incidence angle $\theta=\pi/2$. Contour field of the real part of the normalized velocity potential. Shadow region behind the cylinder is now affected by a strong shoaling effect due to the variable bathymetry}
\label{fig:Cylinder variable depth theta=90}
\end{center}
\end{figure}

\begin{figure}
\begin{center}
\includegraphics[width=8cm]{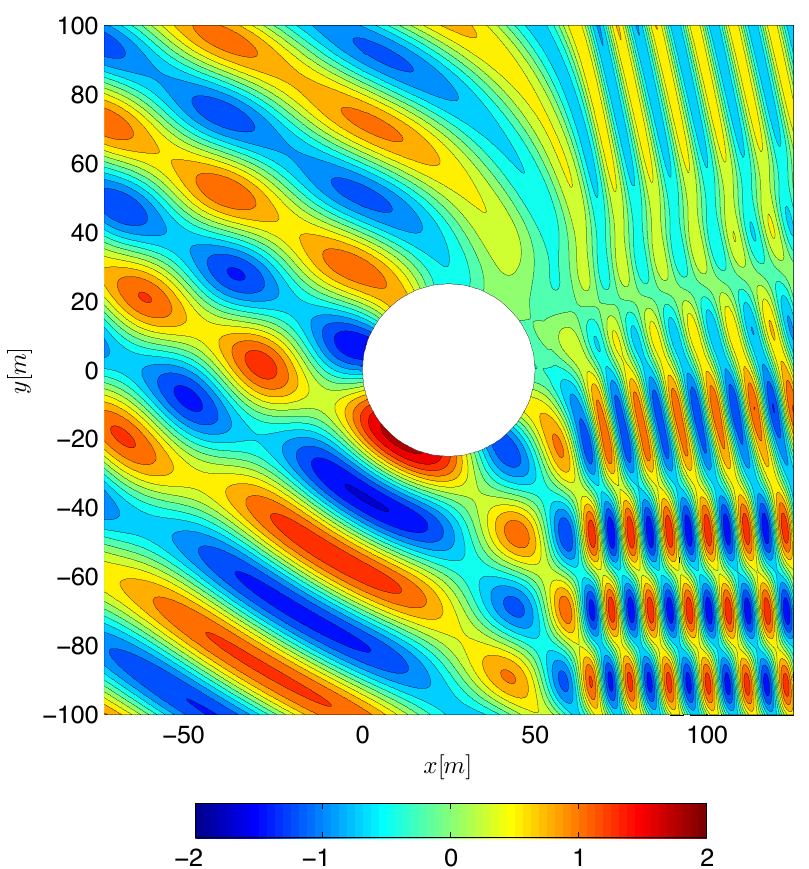}
\caption{BEM solution for a decreasing water depth in the $x$-direction and incidence angle $\theta=\pi/6$. Contour field of the real part of the normalized velocity potential. Combined shoaling and refraction effects behind the cylinder due to the variable bathymetry}
\label{fig:Cylinder variable depth theta=30}
\end{center}
\end{figure}

The problem is next solved for two different incidence angles, $\theta=\pi/2$ and $\theta=\pi/6$, in the same variable water-depth profile considered in Section \ref{SUBSECTION: NM CANAL}. The solution for $\theta=\pi/2$, incident wave normal to the bathymetric lines, is represented in Figure \ref{fig:Cylinder variable depth theta=90}, where a strong shoaling effect is observed behind the cylinder due to the change of water depth in this region. For an inclined incident wave of $\theta=\pi/6$, the solution is shown in Figure \ref{fig:Cylinder variable depth theta=30}. In this case, combined shoaling and refraction effects appear behind the cylinder. Refraction involves a change in the direction of waves as they pass from one medium to another. This phenomenon appears in water waves when traveling from deep to shallow waters, because the waves change their direction and tend to travel perpendicularly to the bathymetric lines. These results demonstrate the ability of the proposed formulation to model the combined effects of shoaling, diffraction and refraction in water wave transmission problems.

\begin{figure}
\begin{center}
\includegraphics[width=0.8\columnwidth]{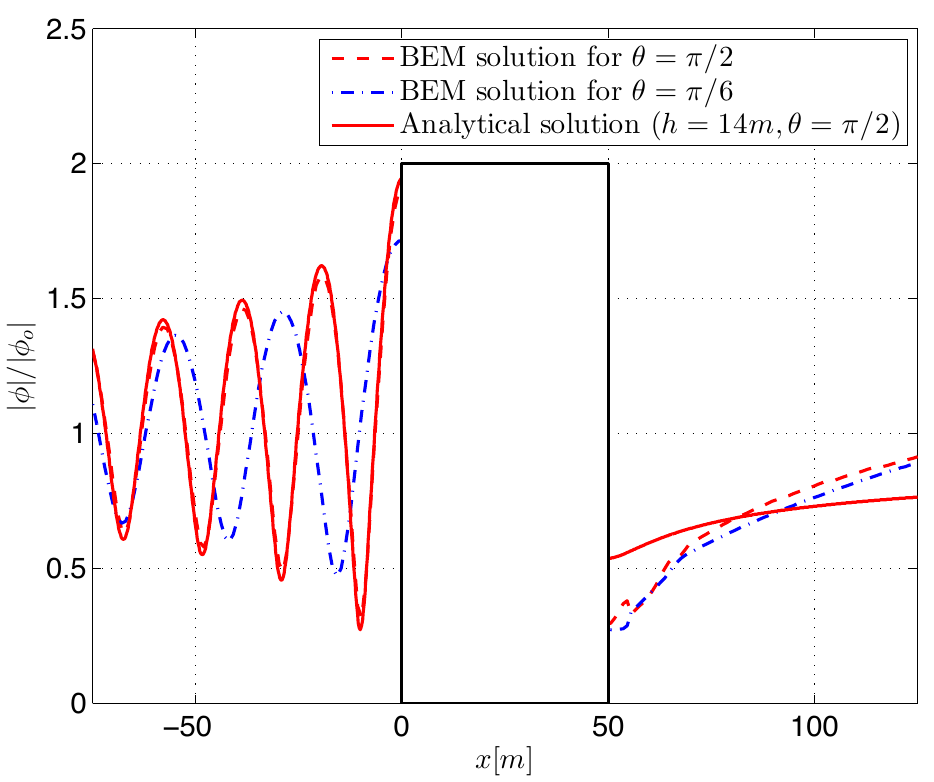}
\caption{Wave scattering by a cylinder in variable water depth. Normalized magnitude of the velocity potential computed with the BEM along the profile located at $y=0$ for two different incidence angles}
\label{fig-Cylinder Variable Bathymetry PROFILE}
\end{center}
\end{figure}

Finally, in Figure \ref{fig-Cylinder Variable Bathymetry PROFILE} we study the influence of variable water-depth and different incidence angles on the normalized velocity potential at $y=0$. We observe that the solution in front of the cylinder is mainly influenced by the incidence angle and that the variable bathymetry is controlling the solution in the wake behind the object. 

\subsection{Elliptic shoal on a sloping bottom}
\label{SUBSECTION: ELLIPTICAL SHOAL}
Although the proposed fundamental solution is restricted to unidirectional variation of the bathymetry, problems with local irregularities of the seabed in two directions can be effectively treated using BEM-FEM coupling techniques \cite{Luis2015} by enclosing the unevenness within a FEM domain that is connected to a BEM model of the external region where the bathymetry varies only in one direction. To demonstrate the effectiveness of this technique and validate the results of the BEM model, we study the scattering produced by an elliptic shoal resting on a sloping seabed. This problem was first studied by Berkhoff et al. \cite{Berkhoff1982} comparing experimental and numerical results based on the MSE and later used by Belibassakis et. al. \cite{Belibassakis2001} to verify a coupled-mode model.

In this problem, the bathymetry is composed of an elliptic shoal superimposed on a sloping bottom with a constant slope of $2\%$. The shape of the inclined background bathymetry $h_i(x,y)$ is given by:
\begin{equation}
h_i(x,y)=\left\lbrace 
\begin{array}{l l}
0.45,& x<-5.85 \\
0.45-0.02(5.85+x), & -5.85 \leq x \leq 14.15 \\
0.05, & x>14.15 
\end{array}
\right.
\end{equation}
and the superimposed shoal, located inside a domain $\Omega_s$ with boundary $\Gamma_s$, produces a disturbance height $h_d(x,y)$ that is evaluated as:
\begin{equation}
h_s(x,y)=0.3-0.5\sqrt{1-\left(\dfrac{x}{3.75}\right)^{2}-\left(\dfrac{y}{5}\right)^{2}}, \quad (x,y) \in \Omega_s
\end{equation}
where the elliptic domain of the disturbance is defined by the condition $\Omega_{s} = \{(x,y) \:| \:(x/3)^{2}+(y/4)^{2} \le 1\}$.

The numerical model consists of an internal region, defined in $\Omega_{s}$, modeled using a FEM approximation of the MSE and an external problem, defined on $\Gamma_s$, modeled by the BEM. The BEM equipped with the proposed fundamental solution allows us to reproduce the infinite domain without the need of using special techniques to satisfy the Sommerfeld radiation condition. A regular mesh of isoparametric quadrilateral finite elements is used to discretize $\Omega_{s}$, with $400$ divisions in the $x$-axis and $300$ divisions in the $y$-axis to have at least $12$ elements per wavelength. The external surface of the domain is meshed using linear two-node boundary elements perfectly matching the finite element discretization on the boundary. Coupling of the FEM and BEM meshes is made node-to-node, compatibilizing nodal velocity-potential and fluxes.

\begin{figure}
\begin{center}
\includegraphics[width=0.85\columnwidth]{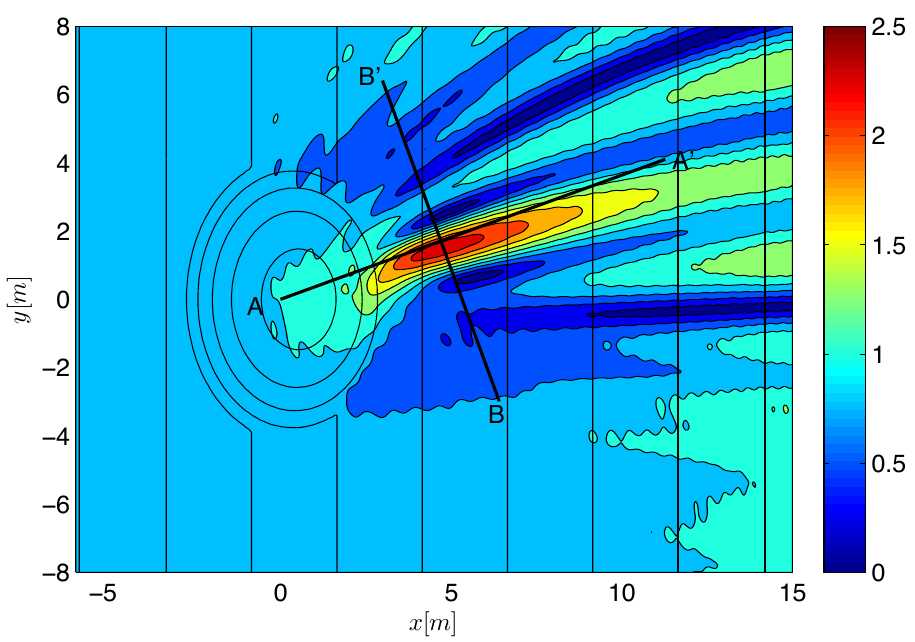} \\
\includegraphics[width=0.85\columnwidth]{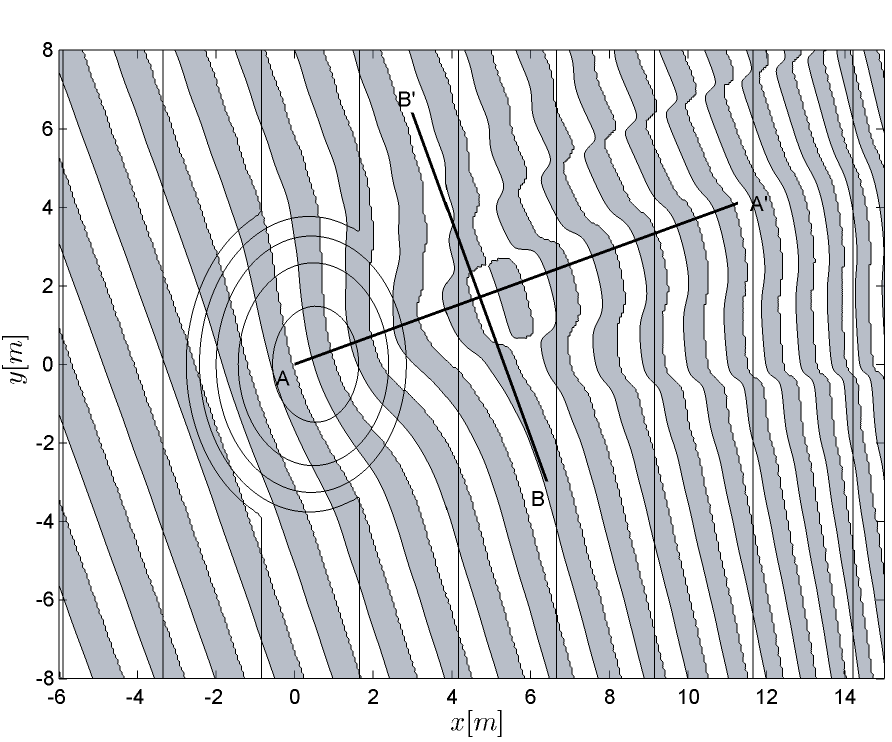}
\caption{Scattering produced by an elliptic bank on a sloping bottom. Absolute value of the normalized wave height (top) and zones of equal phase within an interval of $\pi$ radians (bottom). Two thick lines indicate the location of the sections represented in Figure \ref{fig: sections AB}}
\label{fig: el abs field with sections}
\end{center}
\end{figure}

\begin{figure}
\begin{center}
\includegraphics[width=0.8\columnwidth]{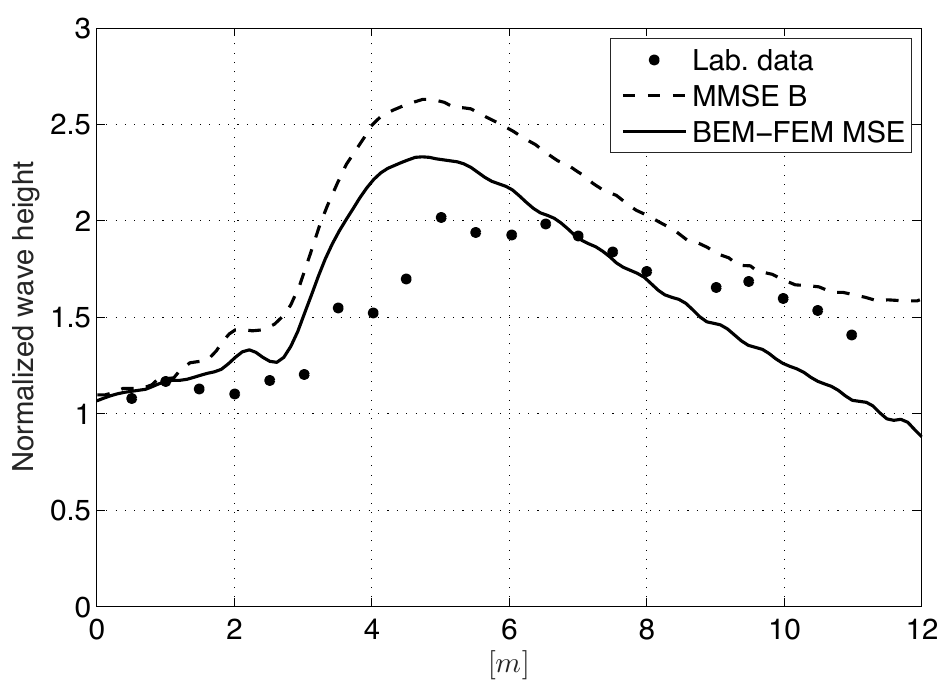} \\
\includegraphics[width=0.8\columnwidth]{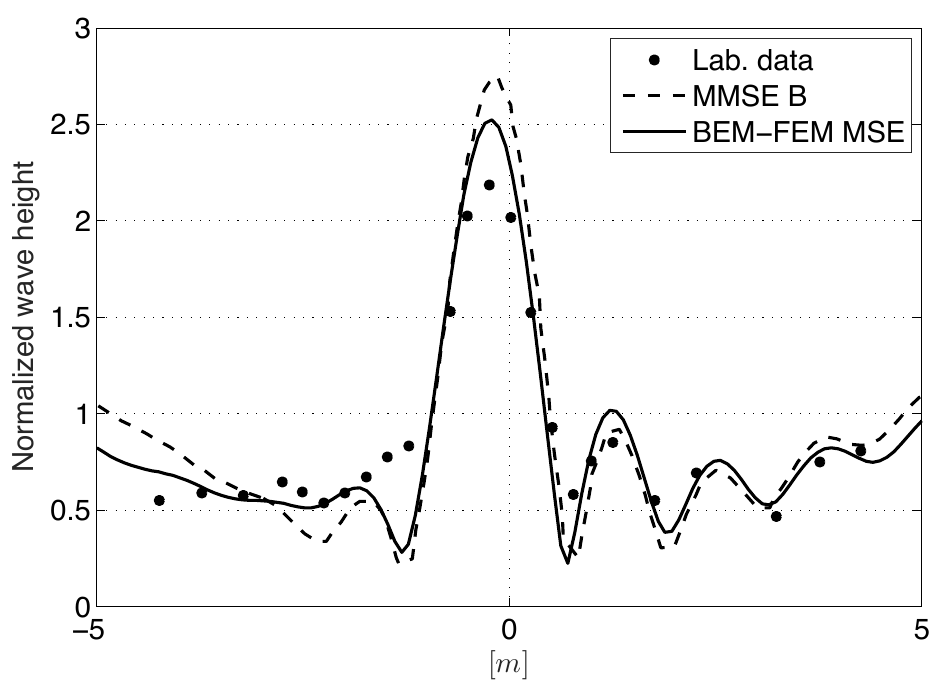}
\caption{Normalized wave height along a longitudinal section A-A' (top) and a transversal section B-B' (bottom), see Figure \ref{fig: el abs field with sections}. Comparison of experimental data from Berkhoff et al. \cite{Berkhoff1982}, numerical results based on the MMSE obtained by Belibassakis et al. \cite{Belibassakis2001} and numerical solution of the MSE obtained by BEM-FEM coupling}
\label{fig: sections AB}
\end{center}
\end{figure}

The considered incoming wave has a period $T=1s$ and enters into the domain with an angle of incidence $\theta=20^{\circ}$ in the $x$ direction. These conditions produce intermediate-water waves with a shallowness ratio varying between $0.30$ and $0.074$ from the deepest area to the shallow region. In Figure \ref{fig: el abs field with sections} (top) we represent the contours of equal-amplitude lines of normalized wave height obtained over the elliptical shoal using the BEM-FEM coupling technique. Two sections are defined in the domain to compare the solution with experimental results obtained by Berkhoff et. al. \cite{Berkhoff1982} and numerical results from Belibassakis et. al \cite{Belibassakis2001} solving the MMSE for this case, see Figure \ref{fig: sections AB}.  Equal phase-range zones of the solution are shown at the bottom of Figure \ref{fig: el abs field with sections}, where it can be appreciated its continuity and absence of spurious reflections at the BEM-FEM interface. 

Sections A-A' and B-B' of the solution are represented in Figure \ref{fig: sections AB}, demonstrating the good behavior of the fundamental solution and presenting a very reasonable agreement with the experimental data. The discrepancies between numerical and experimental results, as the overestimation of the focal peak in the transversal section B-B' and the lack of prediction along the longitudinal section A-A', are mainly attributed to the non-linear effects  \cite{Belibassakis2001}, not considered in the MSE linear theory. The phase results represented in Figure \ref{fig: el abs field with sections} (bottom) are in very good agreement with those measured by Berkhoff et al. \cite{Berkhoff1982}, capturing the presence of two low amplitude points behind the shoal observed in the experiments.

\subsection{Harbor resonance study}
\label{SUBSECTION: NE CHIPIONA HARBOUR}
Finally, in order to explore the possibilities of the proposed BEM formulation in applications with more complex geometries, it is studied the wave transmission problem in a small harbor and its near-by coastal region with variable bathymetry. The harbor is located in the coastal village of Chipiona, at coordinates $36^\circ 44'57''$ latitude and $6^\circ 25'42''$ longitude on the Atlantic shoreline of the C\'{a}diz province, southwest of Spain.

\begin{figure}
\centering
\def\svgwidth{0.8\columnwidth}
\input{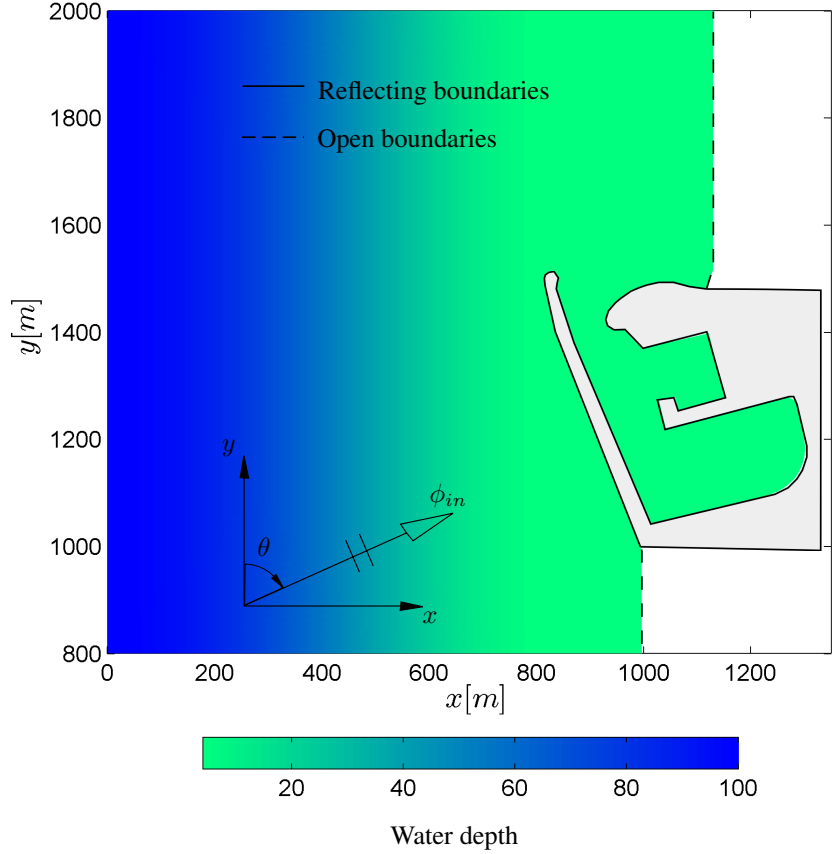}
\caption{Nearshore bathymetric mapping and boundary conditions used for the harbor simulation. Water depth changes along the $x$-coordinate from $100 m$ to $4 m$ between $x=0m$ and $x=840m$. Water depth remains constant for $x \geq 840m$}
\label{fig:HARBOUR_SCHEME}
\end{figure}

The objective of the simulation is to estimate the wave amplification in the interior of the harbor and reproduce the diffraction effects in front of the breakwater. The geometry and boundary conditions of the numerical model are shown in Figure \ref{fig:HARBOUR_SCHEME}, with open boundaries to represent the shoreline and completely reflecting boundaries in the dike and quay walls of the harbor. The open boundaries are modeled assuming a complete absorption of the incoming wave. Bathymetric lines are defined parallel to the shoreline, varying with the polynomial depth function \eqref{eqn:NE:Canal:profile} with coefficients $a_0=100m$, $a_1=0m^{-1}$, $a_2=-4.0816327 \times 10^{-4}m^{-2}$ and $a_3=3.239391 \times 10^{-7}m^{-3}$. In the region, the tidal range is $3.4m$ with a lowest tide of $4m$ that is considered as mean water-depth inside the harbor.

The boundary of the harbor is discretized using 1234 linear boundary elements with at least 20 elements per wave length. As external source, it is considered an incident wave of heigh $1 m$, angle of incidence $\theta=\pi/6$ and a common period value for this area, $T=10s$.

\begin{figure}
\begin{center}
\includegraphics[width=0.8\columnwidth]{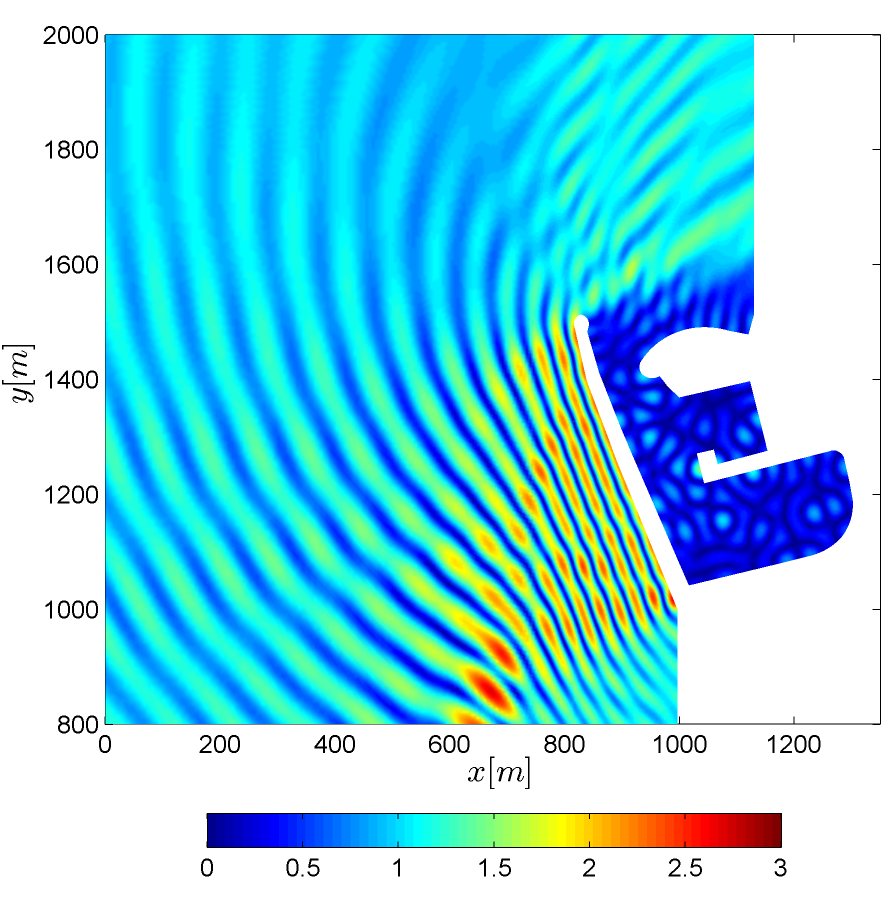}
\caption{WAF diagram calculated using BEM for the harbor of Chipiona example under inclined incident waves and 10-second period wave condition}
\label{fig:HARBOUR WAVE FIELD}
\end{center}
\end{figure}

As part of the post-process and after solving the boundary problem, we computed the solution in the exterior domain at 72539 collocation points that are evenly distributed with a distance between two adjacent points of 1/20 times the local wave length. The WAF obtained for this example is shown in Figure \ref{fig:HARBOUR WAVE FIELD}. It can be observed that some points inside the harbor exhibit a significant amplification factor because dissipation effects, like friction with the sea floor and partial-reflecting boundaries, are not considered in the analysis. However, the use total reflection boundaries is a common practice in harbor resonance studies.

\section{Summary and Conclusions}
\label{SECTION: CONCLUSION}

A complete fundamental solution and BEM formulation for the elliptic Mild-Slope equation in waters of variable depth in one direction has been presented. The Green's function proposed by Belibassakis \cite{Belibassakis2000} has been used as the starting point of our BEM formulation and different examples have been solved to validate the approximation. The main conclusions and findings of this work are the following.

\begin{itemize}
\item The Green's function of Belibassakis \cite{Belibassakis2000} for the MSE with one-directional variable bathymetries has been extended and combined with a boundary element formulation to simulate surface water-wave transmission problems in medium to shallow transition waters.

\item The one-dimensional wave equations in the transformed domain associated with the evaluation of the fundamental solution are solved using a classical Galerkin finite element approximation. It has been observed in our numerical experiments that similar accuracy is obtained using this approach, compared to the second order FD scheme proposed by Belibasakis \cite{Belibassakis2000} for this task.

\item The formulation is able to correctly reproduce the phenomena appearing in water-wave transmission problems: shoaling, diffraction, refraction and the result of their combined effects.

\item Bathymetries with slopes up to 1:3 and contour lines parallel to the shoreline are very common in real problems. The proposed BEM technique allows to simulate these conditions, providing accurate solutions for practical coastal engineering problems.

\item This BEM formulation can be coupled with classical design techniques, like FEM or FDM formulations of the MSE or MMSE \cite{Naserizadeh2011}, to model open sea conditions of variable bathymetry.

\end{itemize}

Finally, it is important to mention that the proposed BEM formulation can also be combined with more advanced partial-reflection boundary conditions and FEM-BEM coupling techniques to extend the range of practical applications. Work in this direction is under way.

\section{Acknowledgements}
This work was supported by the \emph{Ministerio de Econom\'{i}a y Competitividad} of Spain, under the research projects DPI2010-19331and DPI2013-43267-P, which were co-funded by European Regional Development Funds (ERDF).


\providecommand{\bysame}{\leavevmode\hbox to3em{\hrulefill}\thinspace}
\providecommand{\MR}{\relax\ifhmode\unskip\space\fi MR }
\providecommand{\MRhref}[2]{%
  \href{http://www.ams.org/mathscinet-getitem?mr=#1}{#2}
}
\providecommand{\href}[2]{#2}

\end{document}